\pdfoutput=1

\RequirePackage{fix-cm}

\documentclass[%
floatfix,
showkeys,
nofootinbib, %
superscriptaddress, %
]{revtex4-1}

\usepackage{cmap}

\usepackage{ucs}
\usepackage[utf8x]{inputenc}
\usepackage[T1,T2A]{fontenc}
\usepackage[german,russian,english]{babel}

\usepackage[sort&compress]{natbib}
\usepackage{amsmath}
\usepackage{amssymb}

\usepackage{mathtools}
\mathtoolsset{
showonlyrefs,
mathic = true
}

\allowdisplaybreaks

\usepackage{hyperref}
\hypersetup{backref,
 colorlinks=false}
\hypersetup{pdfborder=0 0 0}

\usepackage{microtype}
\UseMicrotypeSet[protrusion]{alltext}

\usepackage{graphicx}

\usepackage[scanall]{psfrag}

\usepackage{listings}
\usepackage{listingsutf8}
\usepackage{fixltx2e}

\usepackage{nicefrac}

\makeatletter
\def\ps@pprintTitle{%
     \let\@oddhead\@empty
     \let\@evenhead\@empty
     \let\@oddfoot\@empty
     \let\@evenfoot\@oddfoot}
\makeatother

\usepackage{mathrsfs}

\begin{document}

\title{Issues in the software implementation of stochastic numerical Runge--Kutta}

\author{Migran N. Gevorkyan}
\email{gevorkyan-mn@rudn.ru}
\affiliation{Department of Applied Probability and Informatics,\\
  Peoples' Friendship University of Russia (RUDN University),\\
  6 Miklukho-Maklaya St, Moscow, 117198, Russian Federation}

\author{Anastasia V. Demidova}
\email{demidova-av@rudn.ru}
\affiliation{Department of Applied Probability and Informatics,\\
  Peoples' Friendship University of Russia (RUDN University),\\
  6 Miklukho-Maklaya St, Moscow, 117198, Russian Federation}

\author{Anna V. Korolkova}
\email{korolkova-av@rudn.ru}
\affiliation{Department of Applied Probability and Informatics,\\
  Peoples' Friendship University of Russia (RUDN University),\\
  6 Miklukho-Maklaya St, Moscow, 117198, Russian Federation}

\author{Dmitry S. Kulyabov}
\email{kulyabov-ds@rudn.ru}
\affiliation{Department of Applied Probability and Informatics,\\
  Peoples' Friendship University of Russia (RUDN University),\\
  6 Miklukho-Maklaya St, Moscow, 117198, Russian Federation}
\affiliation{Laboratory of Information Technologies\\
  Joint Institute for Nuclear Research\\
  6 Joliot-Curie, Dubna, Moscow region, 141980, Russia}

\begin{abstract}
  This paper discusses stochastic numerical methods of Runge-Kutta
  type with weak and strong convergences for systems of stochastic
  differential equations in Itô form. At the beginning we give a brief
  overview of the stochastic numerical methods and information from
  the theory of stochastic differential equations. Then we motivate
  the approach to the implementation of these methods using source
  code generation. We discuss the implementation details and the used
  programming languages and libraries
\end{abstract}


\maketitle

\section{Introduction}
\label{sec:intro}

While modeling technical systems with control it is often required to 
study characteristics of these systems. Also it is necessary to study the influence 
of system parameters on characteristics. In systems with control 
there is a parasitic phenomenon as self-oscillating mode.
We carried out 
studies to determine the region of the self-oscillations emergence. 
However, the parameters of these oscillations were not investigated. 
In this paper, we propose to use the harmonic 
linearization method for this task. This method is used in control 
theory, but this branch of mathematics rarely used in classical 
mathematical modeling. The authors offer a methodological article in order  
to introduce this method to non-specialists.

\section{Introduction}
\label{sec:intro}

The article~\cite{kulyabov:2016:rk-stochastic} describes the Python~\cite{L_Python}
implementation of stochastic numerical Runge-Kutta like methods. This
implementations heavily relies on NumPy and SciPy~\cite{l_scipy}
libraries. We chose Pyhon language because of it's simplicity and
development speed. NumPy's capability to work with multidimensional
arrays as tensors (functions \verb|tensor_dot| and \verb|einsum|) was
also very helpfull. However, the performance was low, and not so much
because of Python slowness, as because we used the large number of
nested loops (up to seven). In this paper, we consider an alternative
approach of stochastic numerical methods implementation, based on
automatical code generation.

This article is divided into three sections. The first section
provides an overview of the main sources and presents information from
the theory of stochastic differential equations (SDE) and methods for
their numerical solution. The second section presents stochastic
numerical schemes for scalar SDE with strong convergence and for SDE
systems with strong and weak convergence. In addition to the general
schemes, several coefficient tables are provided. This allows to
implement a specific numerical method. Finally, the third section
explains the use of code generation for stochastic numerical methods
and describes some details of the generator we have implemented (we
use Jinja2~\cite{l_jinja2} template engine).

\section{Background overview}
\label{sec:review}

In this section, we give a brief overview of the available
publications on stochastic Runge-Kutta methods. We study multistage
numerical schemes without partial derivatives from the drift vector
$\mathbf{f}(t,\mathbf{X})$ and the diffusion matrix
$\mathbf{G}(t,\mathbf{X})$, so we don't consider Milstein
methods~\cite{L_Milstein_1974, L_Milstein_1979, L_Milstein_1986}).

First, who used a stochastic Brownian process for mathematical
modeling was a French mathematician, a student of Henri Poincare ---
Louis J.-B. A. Bachelier (1870--1946) in 1900 in the
work~\cite{L_Bachelier}.

The book by P. E. Kloeden and E. Platen~\cite{L_Kloeden_Platen} is
classical work about numerical methods for SDE. The book provides a
brief introduction to the theory of stochastic Ito and Stratonovich
differential equations and their applications. The last two thirds of
the book are devoted to the presentation of numerical methods in the
sense of strict and weak approximations, including a number of
Runge-Kutta methods.

The dissertation by Andreas Rosler~\cite{L_Andreas_2003} is a
consistent report of stochastic numerical Runge-Kutta-like
methods. The author considers the approximation of Ito and
Stratonovich SDE systems in a weak sense for the scalar and
multidimensional Wiener process. After a brief review of the previous
works, the author develops the stochastic equivalent of labelled trees
theory (labelled trees are used to derive the order conditions in the
case of deterministic Runge-Kutta methods, see, for
example,~\cite{L_Butcher_2003,L_Hairer_2008_en}).

Rossler considers weakly convergent stochastic Runge-Kutta-like
methods for Ito and Stratonovich SDE systems for both the scalar and
the multidimensional Wiener process. In the third and the fifth part
of the dissertation describes specific implementation of the explicit
stochastic numerical methods for weak convergence.

Further results of Rosler studies are described in
articles~\cite{L_Debrabant_2007, L_Debrabant_2013} in collaboration
with K. Debrabant. In the preprint~\cite{L_Debrabant_2013} authors
continue classification of stochastic methods, Runge-Kutta method with
a weak convergence. Several concrete realizations and results of
numerical experiments are given. In the another
preprint~\cite{L_Rossler_2010}, they give tables for fourth stage and
strong order convergence methods $p=3.0$.

Euler--Maruyama method described by Maruyama in the
paper~\cite{L_Maruyama_1955} can be considered as first stochastic
Runge-Kutta-like method. The first systematic study of stochastic
numerical Runge--Kutta-like methods of strong order of convergence
$p=1.0$ is given by V. Rumelin~\cite{L_Rumelin_1982} and E. Platen in
his thesis~\cite{L_Platen_1984}.

Great contribution was made By P. M. Burrage and K. Burrage in a
series of
articles~\cite{L_Burrage_1996,L_Burrage_1997,L_Burrage_1998,L_Burrage_1999,l_burrage_2000}. In
these papers, they not only studied methods of strong order $p=1.5$,
but also extended the theory of labeled trees to the stochastic case.

The article R. Soheili and M. Namjoo~\cite{L_Soheili_2007} obtained
the three methods with strong convergence $p=1.0$ and the numerical
comparison with the method from the book~\cite{L_Kloeden_Platen}.

Some of the first methods with weak convergence are given in the
book~\cite{L_Kloeden_Platen}. Further development they received in
article by Komori and T. Mitsui~\cite{L_Komori_Mitsui_1995} and
in~\cite{L_Mackevicius_1994}. In the article~\cite{L_Tocino_2001} two
three-stage methods, the weak convergence of the $p=2.0$, as well as
numerical experiments are introduced.

In view of the extreme complexity of further improving the order of
accuracy of stochastic numerical schemes, modern studies are devoted
to obtaining numerical schemes for special SDE cases. It is possible
to point out some of such studies about stochastic symplectic
Runge--Kutta-like methods~\cite{L_Burrage2012,L_Ma2015,L_Zhou2017} and
stochastic analogues of the Rosenbrock method~\cite{L_Amiri2017}.

\section{Stochastic Wiener process and software generation of its trajectories}
\label{sec:wiener}

The stochastic process $W (t)$, $t \geqslant 0$ is called scalar
\emph{Wiener process} if the following conditions are
true~\cite{L_Oksendal_en,L_Kloeden_Platen}:
\begin{itemize}
\item $\mathrm{P}\{W(0)=0\} = 1$, or in other words, $W (0) = 0$ is
  almost certain;
\item $W (t)$ is process with independent increments, i.e.
  $\{\Delta W_{i}\}^{N-1}_{0}$ are independent random variables:
  $\Delta W_{I} = W(t_{I+1}) - W (t_{I})$ and
  $0 \leqslant t_{0} < t_{1} < t_{2} < \ldots < t_{N} \leqslant T$;
\item
  $\Delta W_{i} = W (t_{I+1}) - W (t_{I}) \sim \mathcal{N} (0, t_{I +
    1}-t_{I})$ where $0\leqslant t_{I+1} < t_{I} < t$,
  $I=0,1,\ldots, N-1$
\end{itemize}

The symbol $\Delta W_{i} \sim \mathcal{N}(0,\Delta t_{i})$ denotes
that $\Delta W_{i}$ is normally distributed random variable with
expected value $\mathbb{E}[\Delta W_{i}] = \mu = 0$ and variance
$\mathbb{D}[\Delta W_{i}] = \sigma^{2} = \Delta t_{i}$.

The Wiener process is a model of \emph{Brownian motion} (random
walk). If we consider the process $W(t)$ in time points
$0 = t_{0} < t_{1} < t_{2} < \ldots < t_{N-1} < t_{N}$ when it
experiences random additive changes, then directly from the definition
of Wiener process follows:
\[
    W(t_{1}) = W(t_{0}) + \Delta W_0, W(t_{2}) = W(t_{1}) + \Delta W_1, \ldots, W(t_{N}) = W(t_{N-1}) + \Delta W_{N-1},
\]
where $\Delta W_{i} \sim \mathcal{N}(0,\Delta t_{i})$, $\forall i = 0,\ldots,N-1$.

If we write out $W(t_{N})$ as cumulative sum of the increments:
\[
  W(t_{n}) = W(t_{0}) + \sum\limits_{i=0}^{n}\Delta W_{i},\;\;n=0,\ldots,N-1,
\]
and $\mathbb{E}[\Delta W_{I}] = 0$ and
$\mathbb{D}[\Delta W_{i}] = \Delta t_{I}$, then we can show that the
sum of normally distributed random numbers $\Delta W_{I}$ is also a
normally distributed random number:
\[
  \mathbb{E}\sum\limits_{i=0}^{N}\Delta W_{i} = 0, \;\; \mathbb{D}\sum\limits_{i=0}^{n}\Delta W_{i} = \sum\limits_{i=0}^{n} \Delta t_{i} = t_{n} - t_{0},\;\; \sum\limits_{i=0}^{n}\Delta W_{i} \sim \mathcal{N}(0,t_{n}-t_{0}).
\]

The multidimensional Wiener process
$\mathbf{W} (t)\colon \Omega\times[t_{0},T]\to \mathbb{R}^{m}$ is
defined as a random process composed of jointly independent
one-dimensional Wiener processes $W^{1} (t),\ldots,W^{m}
(t)$. Increments of
$\Delta W^{\alpha}_{I},\;\forall \alpha = 1,\ldots,m$ are jointly
independent normally distributed random variables. On the other hand,
the vector $\Delta W^{\alpha}_{I}$ can be represented as a
multidimensional normally distributed random variable with the
expectation vector $\mu = \mathbf{0}$ and the diagonal covariance
matrix.

In the case of a multidimensional stochastic process one has to
generate $m$ sequences of $n$ normally distributed random variables
should be generated.

\section{Stochastic integrals and SDE for scalar Wiener process}

The definition of the stochastic differential equation begins with the
definition of the stochastic integrals.

Let $g(t,x(t))$ be a continuous scalar function depending on the
random process $x(t)$, $W(t)$ is Wiener process, $t\in[t_{0},T]$,
$0 \leqslant t_{0} < t_{1} < t_{2} < \ldots < t_{n-1} < t_{n}
\leqslant T < \infty$ are time points on segment $[t_{0}, T]$, then
\[
I(\theta) = \int\limits_{t_{0}}^{T}g(t,x(t))\mathrm{d} W(t) = \mathop{\mathrm{l.i.m.}\;}\limits_{n\to\infty}\sum^{n-1}_{i=0}g(\theta t_{i+1} + (1-\theta)t_{i}, \theta x(t_{i+1}) + (1-\theta)x(t_{i}))(W(t_{i+1}) - W(t_{i}))
\]
where $\theta\in [0,1]$. $I(\theta)$ is called \textit{stochastic
  integral}. For a more strict definition of stochastic integrals
for a wide class of functions, see~\cite[Chapter 3]{L_Oksendal_en}.

In physics and applied mathematics, two special cases of stochastic
integrals are used:
\begin{itemize}
\item \textit{Ito integrals} (named in honor the Japanese
  mathematician K. Ito)
  \[
  I = I(0) =
  \mathop{\mathrm{l.i.m.}\;}\limits_{n\to\infty}\sum^{n-1}_{i=0}g(t_{i},
  x(t_{i}))(W(t_{i+1}) - W(t_{i})),
  \]
\item \textit{Stratonovich integrals} (named in honor of the soviet
  physicist R. L. Stratonovich)
  \[
  I_{0.5} = I(0.5) =
  \mathop{\mathrm{l.i.m.}\;}\limits_{n\to\infty}\sum^{n-1}_{i=0}g\left(
    \frac{t_{i+1}+t_{i}}{2},\frac{x(t_{i+1}) +
      x(t_{i})}{2}\right)(W(t_{i+1}) - W(t_{i})).
  \]
\end{itemize}

After defining stochastic integrals, we can introduce the integral
equation for the stochastic process $x(t)$~\cite{L_Oksendal_en}.
\[
  x(t_{k}) = x(t_{0}) + \int\limits_{0}^{t}f (\tau,x (\tau))\mathrm{d} \tau + \underbrace {\int\limits_{0}^{t}g(\tau,x (\tau))\mathrm{d} W.} _ {\text{stochastic integral}}
\]
The above integral equation corresponds to the stochastic differential equation (SDE):
\[
  \mathrm{d} x(t) = f(t,x(t))\mathrm{d} t + g(t, x(t))\mathrm{d} W.
\]
It should be noted that the differential $\mathrm{d} x$ is not a usual
``small'' increment of function $x(t)$, but is a random variable. It can
be treated as symbolic notation of the limit of the iterative
scheme. The Wiener process $W (t)$, whose increments are part of the
SDE, is called the driving process for a given SDE.

\section{Ito SDE for multidimensional Wiener process}

Let us consider the probabilistic space
$(\Omega, \mathcal{a}, \mathbb{P})$, where $\Omega$ is the space of
elementary events, $\mathcal{A}$ is sigma-algebra of subsets of the
space $\Omega$, $\mathbb{P}$ is probabilistic measure. The variable
$t$ from the interval $[t_{0}, T] \in \mathbb{R}^{1}$ has physical
sense of time.

Consider the random process
$\mathbf{x} (t) = (x^{1} (t), \ldots, x^{d} (t))^{T}$, where
$\mathbf{x} (t)$ belongs to the functional space
$\mathrm{L}^{2} (\Omega)$ with the norm $\| \cdot \|$. We assume that
the random process $\mathbf{x} (t)$ is a solution for the Ito
SDE~\cite{L_Oksendal_en,L_Kloeden_Platen} if:
\[
  \mathbf{x}(t) = \mathbf{f}(t,\mathbf{x}(t))\mathrm{d} t + \mathbf{G}(t,\mathbf{x}(t))\mathrm{d} \mathbf{W},
\]
where $\mathbf{W} = (W^{1},\dots,W^{m})^T$ is multidimensional Wiener
process, known as \emph{driving} process for SDE.  The function
$\mathbf{f}\colon[t_{0}, T]\times\mathbb{R}^{d} \to \mathbb{R}^{d}$ is
called \emph{drift vector}, and the matrix-valued function
$\mathbf{G}\colon[t_{0}, T]\times\mathbb{R}^{d}\times\mathbb{R}^{m}
\to \mathbb{R}^{d}\times\mathbb{R}^{m}$ is called \emph{diffusion
  matrix}. In addition,
$\mathbf{f} (t,\mathbf{x} (t)) = (f^{1} (t,\mathbf{x}),\ldots,f^{d}
(t,\mathbf{x})^{t}$, and the $\mathbf{G}$ matrix looks like:
\[
  \mathbf{G} = 
  \begin{bmatrix}
    g^{1}_{1}(t,\mathbf{x}) & g^{1}_{2}(t,\mathbf{x}) & \ldots & g^{1}_{m}(t,\mathbf{x})\\
    g^{2}_{1}(t,\mathbf{x}) & g^{2}_{2}(t,\mathbf{x}) & \ldots & g^{2}_{m}(t,\mathbf{x})\\
    \vdots & \vdots & \ddots & \vdots \\
    g^{d}_{1}(t,\mathbf{x}) & g^{d}_{2}(t,\mathbf{x}) & \ldots & g^{d}_{m}(t,\mathbf{x})
  \end{bmatrix}
\]
The same equation can be rewritten in indexed form
\[
  x^{\alpha}(t) = f^{\alpha}(t,x^{\gamma}(t))\mathrm{d} t + \sum\limits_{\beta = 1}^{m}g^{\alpha}_{\beta}(t,x^{\gamma}(t))\mathrm{d} W^{\beta},
\]
where $\alpha,\gamma = 1,\ldots,d$, $\beta = 1,\ldots,m$, and $f^{\alpha}(t,x^{\gamma}(t)) = f^{\alpha}(t,x^{1}(t),\ldots,x^{d}(t))$.

On the interval $[t_{0}, T]$, we introduce the grid
$t_{0}<t_{1}<\ldots<t_{N}=T$ with step $h_{n} = t_{n+1} - t_{n}$,
where $n = 0,\ldots,N-1$ and the maximum grid step
$h = \max{\{h_{n-1}\}^{N}_{1}}$. Next, we assume that the grid is
uniform, then $h_{n} = h = \mathrm{const}$. $\mathbf{x}_{n}$ is grid
function, which approximate a stochastic process $\mathbf{x}(t)$, so
$\mathbf{x}_{0} = \mathbf{x}(t_{0})$,
$\mathbf{x}_{n} \approx \mathbf{x}(t_{n})\;\forall n = 1,\ldots,N$.

\section{Calculation and approximation of multiple Ito integrals of special form}

Here we will not go into the general theory of multiple stochastic Ito
integrals, a reader can refer to the book~\cite{L_Kloeden_Platen} for
additional information. Here we consider multiple special integrals,
which are included in the stochastic numerical schemes.

In General, for the construction of numerical schemes with order of
convergence greater than $p=\frac12$, it is necessary to calculate
single, double and triple Ito integrals of the following form:
\[
  I^{\alpha}(t_{n}, t_{n+1}) = I^{\alpha}(h_{n}) = \int\limits_{t_{n}}^{t_{n+1}} \mathrm{d}W^{\alpha}(\tau),
\]
\[
  I^{\alpha\beta}(t_{n}, t_{n+1}) = I^{\alpha\beta}(h_{n}) = \int\limits_{t_{n}}^{t_{n+1}}\int\limits_{t_{n}}^{\tau_{1}}\mathrm{d}W^{\alpha}(\tau_{2})\mathrm{d}W^{\beta}(\tau_{1}),
\]
\[
  I^{\alpha\beta\gamma}(t_{n}, t_{n+1}) = I^{\alpha\beta\gamma}(h_{n}) = \int\limits_{t_{n}}^{t_{n+1}}\int\limits_{t_{n}}^{\tau_{1}}\int\limits_{t_{n}}^{\tau_{2}}\mathrm{d}W^{\alpha}(\tau_{3})\mathrm{d}W^{\beta}(\tau_{2})\mathrm{d}W^{\gamma}(\tau_{1}),
\]
where $\alpha,\beta,\gamma=0\ldots, m$ and
$W^{\alpha},\alpha=1,\ldots, m$ are components of multidimensional
Wiener process. In the case of $\alpha,\beta,\gamma = 0$, the
increment of $\mathrm{d} W^{0}(\tau)$ is assumed to be
$\mathrm{d} \tau$.

The problem is to get analytical formulas for these integrals with
$\Delta W^{I}_{n} = W^{I}(t_{n+1}) - W^{I}(t_{n})$ in them. Despite
its apparent simplicity, this is not achievable for all possible
combinations of indices. Let us consider in the beginning those cases
when it is possible to obtain an analytical expression, and then turn
to those cases when it is necessary to use an approximating formulas.

In the case of a single integral, the problem is trivial and the
analytic expression can be obtained for any index $\alpha$:
\[
  I^{0}(h_{n}) = \Delta t_{n} = h_{n},\;\; I^{\alpha}(h_{n}) = \Delta W^{\alpha}_{n},\;\alpha=1,\ldots,m.
\]

In the case of a double integral $I^{\alpha\beta} (h_{n})$, the exact
formula takes place only at $\alpha=\beta$:
\[
  I^{00}(h_{n}) = \dfrac{1}{2}\Delta t_{n} = \dfrac{1}{2}h^{2}_{n},\;\; I^{\alpha\alpha}(h_{n}) = \dfrac{1}{2}\left( (\Delta W^{\alpha}_{n})^2 - \Delta t_{n} \right),\;\alpha=1,\ldots,m,
\]
in other cases, when $\alpha\neq \beta$ Express $I^{\alpha\beta} (h_{n})$ by increments of $\Delta W^{\alpha}_{n}$ and $\Delta t_{n}$ in the final form is not possible, so we can only use numerical approximation.

For the mixed case $I^{0\alpha}$ and $I^{\alpha0}$ in~\cite{L_Andreas_2003}, simple formulas of the following form are given:
\[
  I^{0\alpha}(h_{n}) = \dfrac{1}{2}h_{n}\left(I^{\alpha}(h_{n}) - \frac{1}{\sqrt{3}}\zeta^{\alpha}(h_{n})\right),
\]
\[
  I^{\alpha0}(h_{n}) = \dfrac{1}{2}h_{n}\left(I^{\alpha}(h_{n}) + \frac{1}{\sqrt{3}}\zeta^{\alpha}(h_{n})\right),
\]
where $\zeta^{\alpha}_{n}\sim \mathcal{N}(0,h_{n})$ are multidimensional normal distributed random variables.

For the General case $\alpha,\beta = 1,\ldots, m$, the book~\cite{L_Kloeden_Platen} provides the following formulas for approximating the double Ito integral $I^{\alpha\beta}$:
\begin{gather*}
  I^{\alpha\beta}(h_{n}) = \dfrac{\Delta W^{\alpha}_{n} \Delta W^{\beta}_{n} - h_{n}\delta^{\alpha\beta}}{2} + A^{\alpha\beta}(h_{n}),\\
  A^{\alpha\beta}(h_{n}) = \dfrac{h}{2\pi}\sum\limits^{\infty}_{k=1}\dfrac{1}{k}\left[V^{\alpha}_{k}\left(U^{\beta}_{k} + \sqrt{\frac{2}{h_{n}}} \Delta W^{\beta}_{n}\right) - V^{\beta}_{k}\left(U^{\alpha}_{k} + \sqrt{\frac{2}{h_{n}}}\Delta W^{\alpha}_{n}\right)\right],
\end{gather*}
where $V^{\alpha}_{k} \sim \mathcal{N}(0,1)$,
$U^{\alpha}_{k} \sim \mathcal{N}(0,1)$,
$\alpha=1,\ldots,m;\; k=1,\ldots,\infty$; $n=1,\ldots,N$ is numerical
schema number. From the formulas it is seen that in the case
$\alpha=\beta$, we get the final expression for the $I^{\alpha\beta}$,
which we mentioned above. In the case of $\alpha\neq\beta$, one has to
sum the infinite series $a^{\alpha \beta}$. This algorithm gives an
approximation error of order $O(h^2/n)$, where $n$ is number of left
terms of an infinite series $a^{ij}$.

In the article~\cite{l_wiktorsson_2001} a matrix form of approximating
formulas is introduced. Let $\mathbf{1}_{m\times m}$,
$\mathbf{0}_{m\times m}$ be the unit and zero matrices $m\times m$,
then
\begin{gather*}
  \mathbf{I}(h_{n}) = \dfrac{\Delta \mathbf{W}_{n} \Delta \mathbf{W}^{T}_{n} - h_{n}\mathbf{1}_{m\times m}}{2} + \mathbf{A}(h_{n}),\\
  \mathbf{A}(h_{n}) = \dfrac{h}{2\pi}\sum\limits^{\infty}_{k=1}\dfrac{1}{k}\left(\mathbf{V}_{k}(\mathbf{U}_{k} + \sqrt{2/h_{n}} \Delta \mathbf{W}_{n})^{T} - (\mathbf{U}_{k} + \sqrt{2/h_{n}}\Delta \mathbf{W}_{n})\mathbf{V}^{T}_{k}\right),
\end{gather*}
where $\Delta \mathbf{W}_{n}, \mathbf{V}_{k}, \mathbf{U}_{k}$ are independent normally distributed multidimensional random variables:
\begin{gather*}
\Delta \mathbf{W}_{n} = (\Delta W^{1}_{n}, \Delta W^{2}_{n}, \ldots, \Delta W^{m}_{n})^{T} \sim \mathcal{N}(\mathbf{0}_{m\times m}, h_{n}\mathbf{1}_{m\times m}),\\
\mathbf{V}_{k} = (V^{1}_{k}, V^{2}_{k},\ldots, V^{m}_{k})^{T} \sim \mathcal{N}(\mathbf{0}_{m\times m}, \mathbf{1}_{m\times m}),\;\;
\mathbf{U}_{k} = (U^{1}_{k}, U^{2}_{k},\ldots, U^{m}_{k})^{T} \sim \mathcal{N}(\mathbf{0}_{m\times m}, \mathbf{1}_{m\times m}).
\end{gather*}
If the programming language supports vectored operations with multidimensional arrays, these formulas can provide a benefit to the performance of the program.

Finally, consider a triple integral. In the only numerical scheme in
which it occurs, it is necessary to be able to calculate only the case
of identical indexes $\alpha=\beta=\gamma$. For this case,
~\cite{L_Andreas_2003} gives the following formula:
\[
  I^{\alpha\alpha\alpha}(h_{n}) = \dfrac{1}{6}\left((I^{\alpha}(h_{n}))^{3} - 3 I^{0}(h_{n})I^{\alpha}(h_{n})\right) = \dfrac{1}{6}\left((\Delta W^{\alpha}_{n})^{3} - 3 h_{n}\Delta W^{\alpha}_{n}\right).
\]

\section{Strong and weak convergence of the approximating function}

Before proceeding to the formulation of numerical schemes, it is
necessary to determine the criterion of accuracy of approximation of
the simulated process $\mathbf{x}(t)$ by the grid function
$\mathbf{x}_{n}$. Two criteria are used: \emph{weak} and \emph{strong}
convergence.

The sequence of approximating functions $\{\mathbf{x}_{n}\}^{N} _ {1}$
converges with order $p$ to the exact solution $\mathbf{x}(t)$ of SDE
in moment $T$ in \emph{strong sense} if constant $C > 0$ exists and
$\delta_{0} > 0$ such as $\forall h \in (0,\delta_{0}$ and following
condition is fulfilled:
\[
\mathbb{e}(\|\mathbf{x}(T) - \mathbf{x}_{N}\|) \leqslant Ch^{p}.
\]

The sequence of approximating functions $\{\mathbf{x}_{n}\}^{N} _ {1}$
converges with order $p$ to the exact solution $\mathbf{x}(t)$ of SDE
in moment $T$ in \emph{weak sense} if constant $C_{F}>0$ exists and
$\delta_{0}>0$ such as $\forall h \in (0,\delta_{0}]$ and the
following condition is fulfilled:
\[
\left|\mathbb{E}\left[F(\mathbf{x}(T))\right] - \mathbb{E}\left[F(\mathbf{x}_{N})\right]\right| \leqslant C_{F}h^{p}.
\]

Here $F \in C^{2(p+1)}_{\mathrm{P}}(\mathbb{R},\mathbb{R}^{d})$ is a
continuous differentiable functional with polynomial growth.

If the $\mathbf{G}$ matrix is zero, then the strong convergence
condition is equivalent to the deterministic case, but the order of
strong convergence is not necessarily a natural number and can take
fractional-rational values.

It is important to note that the choice of the convergence type
depends on the problem one has to solve. Increasing the order of
strict convergence leads to more accurate approximation of the
trajectories of $\mathbf{x}(t)$. If one wants to calculate, for
example, the moment of a random process $\mathbf{x}(t)$ or a
generalized functional of the form $\mathbb{E}[F(\mathbf{x}(t))]$, one
should increase the order of weak convergence.

\section{Stochastic Runge--Kutta-like numerical methods}
\subsection{Euler--Maruyama numerical method}

The simplest numerical method for solving scalar equations and systems
of SDEs is the Euler--Maruyama method, named in honor of Gisiro
Maruyama, which extended the classical Euler method for ODEs to the
case of equation~\cite{L_Maruyama_1955}. The method is easily
generalized to the case of multidimensional Wiener process.
\begin{align*}
  &x^{\alpha}_{0} = x^{\alpha}(t_{0}),\\
  &x^{\alpha}_{n+1} = x^{\alpha}_{n} + f^{\alpha}(t_{n}, x^{\alpha}_{n})h_{n} + \sum\limits^{d}_{\gamma=1}G^{\alpha}_{\beta}(t_{n},x^{\gamma}_{n})\Delta W^{\beta}_{n}.
\end{align*}
From the formula we can see, that each step requires only
corresponding to this step increment $\Delta W^{\beta}_{n}$. The
method has a strong order $(p_{d},p_{s}) = (1.0, 0.5)$. The value
$p_{d}$ denotes the deterministic accuracy order, when the method is
used for the equation with $G(t,x^{\alpha}(t))\equiv 0$. The value
$p_{s}$ denotes the stochastic part approximation order.

\subsection{Weak stochastic Runge--Kutta-like method with order $1.5$ for a scalar Wiener process}
In the case of a scalar SDE, the drift vector
$f^{\alpha}(t, x^{\gamma})$ and the diffusion matrix
$G^{\alpha}_{\beta}(t, x^{\gamma})$ become $f(t, x)$ and $g(t, x)$
scalar functions, and the driving Wiener process $W^{\beta}_t$ is
scalar $W_t$. For scalar SDE it is possible to construct a numerical
scheme with strong convergence $p=1.5$:
\begin{align*}
  &X^{i}_{0} = x_{n} + \sum\limits^{s}_{j=1}A_{0j}^{i}f(t_{n} + c^{j}_{0}h_{n}, X^{j}_{0})h_{n} + \sum\limits^{s}_{j=1}B^{i}_{0j}g(t_{n} + c^{j}_{1}h_{n}, X^{j}_{1})\dfrac{I^{10}(h_{n})}{\sqrt{h_{n}}},\\
  &X^{i}_{1} = x_{n} + \sum\limits^{s}_{j=1}A_{1j}^{i}f(t_{n} + c^{j}_{0}h_{n}, X^{j}_{0})h_{n} + \sum\limits^{s}_{j=1}B^{i}_{1j}g(t_{n} + c^{j}_{1}h_{n}, X^{j}_{1})\sqrt{h_{n}},\\
  &x_{n+1} = x_{n} + \sum\limits^{s}_{i=1}a_{i}f(t_{n}+c^{i}_{0}h_{n}, X^{i}_{0})h_{n} + \\
  & \qquad + \sum\limits^{s}_{i=1}\left(b^{1}_{i}I^{1}(h_{n}) + b^{2}_{i}\dfrac{I^{11}(h_{n})}{\sqrt{h_{n}}} + b^{3}_{i}\dfrac{I^{10}(h_{n})}{h_{n}} + b^{4}_{i}\dfrac{I^{111}(h_{n})}{h_{n}}\right)g(t_{n} + c^{i}_{1}h_{n},X^{i}_{1}),
\end{align*}
where $i,j=1,\ldots,s$ ($s$ is numbers of method's stages). The generalized Butcher table~\cite{L_Rossler_2010} has fallowing form:
\[
  {\renewcommand{\arraystretch}{1.5}%
  \begin{array}{c|c|c|c}
    c^{i}_{0} & A^{i}_{0j} & B^{i}_{0j} &\\ \hline
    c^{i}_{1} & A^{i}_{1j} & B^{i}_{1j} &\\ \hline
    & a_{i} & b^{1}_{i} & b^{2}_{i}\\ \hline
    &  & b^{3}_{i} & b^{4}_{i}\\
  \end{array}}
\]

In the above numerical scheme, the Wiener stochastic process is
present in implicit way. It is "hidden" inside the stochastic Ito
integrals: $I^{10}(h_{n})$, $I^{1}(h_{n})$, $I^{11}(h_{n})$,
$I^{111}(h_{n})$. For scalar case they are simplified:
\begin{align*}
  & I^{1}(h_{n}) = \Delta W_{n},\\
  & I^{10}(h_{n}) = \dfrac{1}{2}h_{n}(\Delta W_{n} + \zeta_{n}/\sqrt{3}),\\
  & I^{11}(h_{n})  = \dfrac{1}{2}((\Delta W_{n})^2 - h_{n}), \\
  & I^{111}(h_{n}) = \dfrac{1}{6}((\Delta W_{n})^3 - h_{n}\Delta W_{n}),
\end{align*}
where $\zeta_{n}\sim \mathcal{N}(0,h_{n}), \sigma = \sqrt{h_{n}}$.

Rossler introduce two Butcher tables for strong scalar methods in
preprint~\cite{L_Rossler_2010} for $s = 4$
\[
  \text{\texttt{SRK1W1}: }
  {\renewcommand{\arraystretch}{1.1}%
\begin{array}{c|cccc|cccc|cccc}
  0 & 0 & 0 & 0 & 0 & 0 & 0 & 0 & 0 &  &  &  & \\
  3/4 & 3/4 & 0 & 0 & 0 & 3/2 & 0 & 0 & 0 &  &  &  & \\
  0 & 0 & 0 & 0 & 0 & 0 & 0 & 0 & 0 &  &  &  & \\
  0 & 0 & 0 & 0 & 0 & 0 & 0 & 0 & 0 &  &  &  & \\
  \hline
  0 & 0 & 0 & 0 & 0 & 0 & 0 & 0 & 0 &  &  &  & \\
  1/4 & 1/4 & 0 & 0 & 0 & 1/2 & 0 & 0 & 0 &  &  &  & \\
  1 & 1 & 0 & 0 & 0 & -1 & 0 & 0 & 0 &  &  &  & \\
  1/4 & 0 & 0 & 1/4 & 0 & -5 & 3 & 1/2 & 0 &  &  &  & \\
  \hline
  & 1/3 & 2/3 & 0 & 0 & -1 & 4/3 & 2/3 & 0 & -1 & 4/3 & -1/3 & 0\\
  \hline
  &   &   &   &   & 2 & -4/3 & -2/3 & 0 & -2 & 5/3 & -2/3 & 1
\end{array}}
\]
\[
  \text{\texttt{SRK2W1}: }
  {\renewcommand{\arraystretch}{1.1}%
  \begin{array}{c|cccc|cccc|cccc}
  0 & 0 & 0 & 0 & 0 & 0 & 0 & 0 & 0 &  &  &  & \\
  1 & 1 & 0 & 0 & 0 & 0 & 0 & 0 & 0 &  &  &  & \\
  1/2 & 1/4 & 1/4 & 0 & 0 & 1 & 1/2 & 0 & 0 &  &  &  & \\
  0 & 0 & 0 & 0 & 0 & 0 & 0 & 0 & 0 &  &  &  & \\
  \hline
  0 & 0 & 0 & 0 & 0 & 0 & 0 & 0 & 0 &  &  &  & \\
  1/4 & 1/4 & 0 & 0 & 0 & -1/2 & 0 & 0 & 0 &  &  &  & \\
  1 & 1 & 0 & 0 & 0 & 1 & 0 & 0 & 0 &  &  &  & \\
  1/4 & 0 & 0 & 1/4 & 0 & 2 & -1 & 1/2 & 0 &  &  &  & \\
  \hline
  & 1/6 & 1/6 & 2/3 & 0 & -1 & 4/3 & 2/3 & 0 & -1 & -4/3 & 1/3 & 0\\
  \hline
  &   &   &   &   & 2 & -4/3 & -2/3 & 0 & -2 & 5/3 & -2/3 & 1
\end{array}}
\]
The numerical schema for the first table we denote as \texttt{SRK1W1}
and the second as \texttt{SRK2W2}. Methods \texttt{SRK1W1} and
\texttt{SRK2W1} have strong orders $(p_{d},p_{s}) = (2.0, 1.5)$ and
$(p_{d},p_{s}) = (3.0, 1.5)$ respectively. One more method $p_{s}=1.0$
is introduced in book~\cite{L_Kloeden_Platen} and its Butcher table
has the following form:
\[
  \text{\texttt{KlPl}: }
  {\renewcommand{\arraystretch}{1.1}%
  \begin{array}{c|cc|cc|cc}
  0&0&0&0&0&&\\
  0&0&0&0&0&&\\ \hline
  0&0&0&0&0&&\\
  0&1&0&1&0&&\\ \hline
  0&1&0&1&0&-1&1\\\hline
   & & &1&0&0&0\\
  \end{array}}
\]

\subsection{Stochastic Runge--Kutta method with strong order $p=1.0$ for vector Wiener process}
For SDE system with a multidimensional Wiener process, one can
construct a stochastic numerical Runge-Kutta scheme of strong order
$p_{s} = 1.0$ using single and double Ito
integrals~\cite{L_Rossler_2010}.
\begin{align*}
  &X^{0i\alpha} = x^{\alpha}_{n} + \sum\limits^{s}_{j=1}A^{i}_{0j}f^{\alpha}(t_{n} + c^{j}_{0}h_{n}, X^{0j\beta})h_{n} + \sum\limits_{l=1}^{m}\sum\limits^{s}_{j=1}B^{i}_{0j}G^{\alpha}_{l}(t_{n} + c^{j}_{1}h_{n},X^{lj\beta})I^{l}(h_{n}),\\
  &X^{ki\alpha} = x^{\alpha}_{n} + \sum\limits^{s}_{j=1}A^{i}_{1j}f^{\alpha}(t_{n} + c^{j}_{0}h_{n},X^{0j\beta})h_{n} + \sum\limits_{l=1}^{m}\sum\limits^{s}_{j=1}B^{i}_{1j}G^{\alpha}_{l}(t_{n} + c^{j}_{1}h_{n},X^{lj\beta})\dfrac{I^{lk}(h_{n})}{\sqrt{h_{n}}},\\
  &x^{\alpha}_{n+1} = x^{\alpha}_{n} + \sum\limits^{s}_{i=1}a_{i}f^{\alpha}(t_{n} + c^{i}_{0}h_{n},X^{0i\beta})h_{n} + \sum\limits^{m}_{k=1}\sum\limits^{s}_{i=1}(b^1_{i}I^{k}(h_{n}) + b^2_{i}\sqrt{h_{n}})G^{\alpha}_{k}(t_{n} + c^{i}_{1}h_{n},X^{ki\beta}),
\end{align*}
$n=0,1,\ldots,N-1$; $i=1,\ldots,s$; $\beta,k=1,\ldots,m$; $\alpha = 1,\ldots,d$. Its generalized Butcher table has the following form~\cite{L_Rossler_2010}:
\[
  {\renewcommand{\arraystretch}{1.4}%
  \begin{array}{c|c|c|c}
    c^{i}_0 & A^{i}_{0j} & B^{i}_{0j} &\\ \hline
    c^{i}_1 & A^{i}_{1j} & B^{i}_{1j} &\\ \hline
    & a_{i} & b^{1}_{i} & b^{2}_{i} \\
  \end{array}}
\]

Rossler introduce two Butcher tables for strong scalar methods in
preprint~\cite{L_Rossler_2010} for $s = 3$
\[
  \text{\texttt{SRK1Wm}: }
  {\renewcommand{\arraystretch}{1.1}%
  \begin{array}{c|ccc|ccc|ccc}
    0 & 0 & 0 & 0 & 0 & 0 & 0 &  &  & \\
    0 & 0 & 0 & 0 & 0 & 0 & 0 &  &  & \\
    0 & 0 & 0 & 0 & 0 & 0 & 0 &  &  & \\
    \hline
    0 & 0 & 0 & 0 & 0 & 0 & 0 &  &  & \\
    0 & 0 & 0 & 0 & 1 & 0 & 0 &  &  & \\
    0 & 0 & 0 & 0 & -1 & 0 & 0 &  &  & \\
    \hline
    & 1 & 0 & 0 & 1 & 0 & 0 & 0 & 1/2 & -1/2\\
  \end{array}}
\]
\[
  \text{\texttt{SRK2Wm}: }
  {\renewcommand{\arraystretch}{1.1}%
    \begin{array}{c|ccc|ccc|ccc}
    0 & 0 & 0 & 0 & 0 & 0 & 0 &  &  & \\
    1 & 1 & 0 & 0 & 0 & 0 & 0 &  &  & \\
    0 & 0 & 0 & 0 & 0 & 0 & 0 &  &  & \\
    \hline
    0 & 0 & 0 & 0 & 0 & 0 & 0 &  &  & \\
    1 & 1 & 0 & 0 & 1 & 0 & 0 &  &  & \\
    1 & 1 & 0 & 0 & -1 & 0 & 0 &  &  & \\
    \hline
    & 1/2 & 1/2 & 0 & 1 & 0 & 0 & 0 & 1/2 & -1/2\\
  \end{array}}
\]
Methods \texttt{SRK1Wm} and \texttt{SRK2Wm} have strong order $(p_{d},p_{s}) = (1.0, 1.0)$ and $(p_{d},p_{s}) = (2.0, 1.0)$.

\subsection{Stochastic Runge--Kutta method with weak order $p=2.0$ for vector Wiener process}

Numerical methods with weak convergence are good for approximation the
distribution characteristics of stochastic process
$x^{\alpha}(t)$. The weak numerical method does not need information
about the trajectory of driving Wiener process $W^{\alpha}_{n}$ and
random increments for these methods can be generated on another
probability space.
\begin{align*}
& X^{0i\alpha} = x^{\alpha}_{n} + \sum\limits^{s}_{j=1}A_{0j}^{i}f^{\alpha}(t_{n} + c^{j}_{0}h_{n}, X^{0j\beta})h_{n} + \sum\limits^{s}_{j=1}\sum\limits^{m}_{l=1}B_{0j}^{i}G^{\alpha}_{l}(t_{n}+c^{j}_{1}h_{n}, X^{lj\beta})\hat{I}^{l},\\
& X^{ki\alpha} = x^{\alpha}_{n} + \sum\limits^{s}_{j=1}A^{i}_{1j}f^{\alpha}(t_{n}+c^{j}_{0}h_{n}, X^{0j\beta})h_{n} + \sum\limits^{s}_{j=1}B^{i}_{1j}G^{\alpha}_{k}(t_{n} + c_{1}^{j}h_{n}, X^{kj\beta})\sqrt{h_{n}},\\
& \widehat{X}^{ki\alpha} = x^{\alpha}_{n} + \sum\limits^{s}_{j=1}A^{i}_{2j}f^{\alpha}(t_{n}+c^{j}_{0}h_{n}, X^{0j\beta})h_{n} + \sum\limits^{s}_{j=1}\sum\limits^{m}_{l=1,l\neq k}B^{i}_{2j}G^{\alpha}_{l}(t_{n} + c_{1}^{j}h_{n}, X^{lj\beta})\frac{\hat{I}^{kl}}{\sqrt{h_{n}}},\\
& x^{\alpha}_{n+1} = x^{\alpha}_{n} + \sum\limits^{s}_{i=1}a_{i}f^{\alpha}(t_{n}+c^{i}_{0}, X^{0i\beta})h_{n} + \sum\limits^{s}_{i=1}\sum\limits^{m}_{k=1}\left(b^{1}_{i}\hat{I}^{k} + b^{2}_{i}\frac{\hat{I}^{kk}}{\sqrt{h_{n}}}\right) G^{\alpha}_{k}(t_{n}+c^{i}_{1}h_{n}, X^{ki\beta}) + \\
& +  \sum\limits^{s}_{i=1}\sum\limits^{m}_{k=1}\left(b^{3}_{i}\hat{I}^{k} + b^{4}_{i}\sqrt{h_{n}}\right)G^{\alpha}_{k}(t_{n}+c^{i}_{2}h_{n}, \widehat{X}^{ki\beta})
\end{align*}
Generalized Butcher table has the following form~\cite{L_Rossler_2010}
\[
  {\renewcommand{\arraystretch}{1.5}%
  \begin{array}{c|c|c|c}
    c^{i}_{0} & A^{i}_{0j} & B^{i}_{0j} &\\ \hline
    c^{i}_{1} & A^{i}_{1j} & B^{i}_{1j} &\\ \hline
    c^{i}_{2} & A^{i}_{2j} & B^{i}_{2j} &\\ \hline
    & a_{i} & b^{1}_{i} & b^{2}_{i}\\ \hline
    &  & b^{3}_{i} & b^{4}_{i}\\
  \end{array}}
\]

From the paper~\cite{L_Debrabant_2013} we get two Butcher tables:
{
\[
  {\renewcommand{\arraystretch}{1.5}%
  \begin{array}{c|ccc|ccc|ccc}
    0&0&0&0&0&0&0&&&\\
    1&1&0&0&\frac13&0&0&&&\\
    \frac{5}{12}&\frac{25}{144}&\frac{35}{144}&0&-\frac56&0&0&&&\\ \hline
    0&0&0&0&0&0&0&&&\\
    \frac14&\frac14&0&0&\frac12&0&0&&&\\
    \frac14&\frac14&0&0&-\frac12&0&0&&&\\ \hline
    0&0&0&0&0&0&0&&&\\
    0&0&0&0&1&0&0&&&\\
    0&0&0&0&-1&0&0&&&\\ \hline
    &\frac{1}{10}&\frac{3}{14}&\frac{24}{35}&1&-1&-1&0&1&-1\\ \hline
    &&&&\frac12&-\frac14&-\frac14&0&\frac12&-\frac12\\
  \end{array}}\quad
  {\renewcommand{\arraystretch}{1.5}%
  \begin{array}{c|ccc|ccc|ccc}
    0&0&0&0&0&0&0&&&\\
    1&1&0&0&1&0&0&&&\\
    0&0&0&0&0&0&0&&&\\ \hline
    0&0&0&0&0&0&0&&&\\
    1&1&0&0&1&0&0&&&\\
    1&1&0&0&-1&0&0&&&\\ \hline
    0&0&0&0&0&0&0&&&\\
    0&0&0&0&1&0&0&&&\\
    0&0&0&0&-1&0&0&&&\\ \hline
    &\frac{1}{2}&\frac{1}{2}&0&\frac12&\frac14&\frac14&0&\frac12&-\frac12\\ \hline
    &&&&-\frac12&\frac14&\frac14&0&\frac12&-\frac12\\
  \end{array}}
\]}

In the weak numerical schema $\hat{I}^{kl}$ are
\[
  \hat{I}^{kl} = 
  \left\{
    \begin{aligned}
      &\dfrac{1}{2}(\hat{I}^{k}\hat{I}^{l} - \sqrt{h}_{n}\tilde{I}^{k}),\; k<l,\\
      &\dfrac{1}{2}(\hat{I}^{k}\hat{I}^{l} + \sqrt{h}_{n}\tilde{I}^{l}),\; l<k,\\
      &\dfrac{1}{2}((\hat{I}^{k})^2 - h_{n}).\; k=l.
    \end{aligned}
  \right.
\]

Where $\hat{I}^{k}$ denotes three point distributed random
variable. It means, that $\hat{I}^{k}$ may have three values
$ \{-\sqrt{3h_{n}}, 0, \sqrt{3h_{n}}\}$ with probabilities $1/6$,
$2/3$ and $1/6$ respectively. $\tilde{I}^{k}$ denotes two point
distributed random variable $\{-\sqrt{h_{n}}, \sqrt{h_{n}}\}$ with
probabilities $1/2$ and $1/2$ respectively.

\section{Analysis of implementation difficulties of stochastic Runge--Kutta numerical methods}

As can be seen from the formulas, stochastic Runge-Kutta methods are
much more complicated than their classical analogues. In addition to
the cumbersome formulas, we can highlight the following factors that
complicate the implementation stochastic methods in software, as well
as their application to the numerical solution of SDEs.
\begin{itemize}
\item When choosing a particular method, it is necessary to consider
  what type of convergence is necessary to provide for this particular
  problem, as well as which of the stochastic equations should be
  solved --- in Ito or Stratonovich form. This increases the number of
  algorithms one has to implement.
\item For methods with strong convergence of greater then one at each
  step it is necessary to solve the resource-intensive problem of
  stochastic integrals approximation.
\item In the numerical scheme, there are not only matrices and
  vectors, but also tensors (four-dimensional arrays) with which it is
  necessary to perform a convolution operation on several indexes. The
  implementation of convolution via summation using normal cycles
  results in a significant performance drop.
\item Weak methods requires the Monte Carlo and, therefore, a large
  number of repeated computations of the numerical solution. Since the
  Monte Carlo method converges approximately as $1 / \sqrt{N}$, where
  $n$ --- number of calculations, to achieve an accuracy of at least
  $10^{-3}$, it is necessary to perform minimum $10^6$ tests.
\end{itemize}

The most significant performance drop occurs when implementing a
universal algorithm, that is, a program that can make a calculation
using an arbitrary coefficient table. In this case, we have to use a
large number of nested loops in order to organize the summation. The
presence in the schemes of double sums and complex combination of
indices in the multipliers under the sign of these sums complicates
complicates the implementation even more and the number of nested
cycles increases to six.  In addition to these specific features, we
mention a few reasons for the performance drop, which also take place
in case of deterministic numerical methods. The obvious way to store
the coefficients of the methods is to use arrays. However, in explicit
methods that we consider, the matrix is lower-diagonal and storing it
as a two-dimensional array results in more than half of the allocated
memory being spent on storing zeros.

If you examine the source codes of popular routines that implement
classical explicit embedded Runge--Kutta methods, one may find that
these programs use a set of named constants rather than arrays to
store the coefficients of the method. It is also caused by the fact
that the operations with scalar variables in most programming
languages are faster than operations on arrays.

We wish to preserve the requirement of code universality and at the
same time to increase the speed of calculations and reduce the memory
consumption. That led us to automatic code generation from one
template.

In addition to performance gains, automatic code generation allows you
to add or modify all functions at once by editing only one
template. This allows both to reduce the number of errors and to
generate different variants of functions for different purposes.

\section{Automatic code generation}

For code generation we use Python 3 language. The program is open
source and available on bitbucket repository by URL
\url{bitbucket.org/mngev/sde_num_generation}. The repository contains
module \texttt{stochastic}. This module implements Wiener stochastic
process and the numerical methods we considered in this paper. Most
part of the module's code are generated by scripts from
\texttt{generator} directory.

For the the code generation, we used Jinja2~\cite{l_jinja2} template
engine. This library was originally developed to generate HTML pages,
but it has a very flexible syntax and can be used as a universal tool
for generating text files of any kind, including source codes in any
programming languages. In addition to Jinja2, we also used NumPy
library to work with arrays and speed-up some calculations.

In addition to the two external libraries listed above, the standard
\texttt{fraction} module was used. It allows to specify the
coefficients of the method as rational fractions, and then convert
them to float type with the desired order of accuracy. Also we use
\texttt{typing} module to annotate the types of function arguments
(Python 3.5 and above feature).

Templates are files with Python source code with insertions of Jinja2
specific commands. Information about the coefficients of the methods
is stored separately, in a structured form of JSON format. This makes
it easy to add new methods and modify old ones by editing  JSON
files. Currently we use methods with coefficients presented
in~\cite{L_Debrabant_2007,L_Debrabant_2013,L_Kloeden_Platen}.

Python itself is used as the language for already generated functions
with the active use of NumPy library, which allows to get acceptable
performance. However, the generated code can be easily reformatted to
match the syntax of any other programming language. We plan to modify
the program to generate code in Julia language (julialang.org). This
language was introduced in 2012 and initially focused on scientific
computing. Currently, he is intensively developing and gaining
popularity. To date, the current version is 0.6.2. Julia provides
performance comparable to C++ and Fortran, but it is a dynamic
language with interactive command line (REPL) capability similar to
IPython and can be integrated into an interactive Jupyter environment.

The current version of the library exceeds the one described by the
authors in~\cite{kulyabov:2016:rk-stochastic}. The use of
auto-generation made it possible not to use nested loops, which
reduced the number of memory allocations, and greatly simplified the
code.

\subsection{Realisation of automatic code generation}
To study the calculation errors and the efficiency of different
stochastic numerical methods, it is necessary to have a universal
implementation of such methods. The universality means the possibility
to use any stochastic method with a desired strong or weak error by
setting its coefficient table. With direct transfer of mathematical
formulas to the program code, one need to use about five nested
cycles, which extremely reduces performance, since such code does not
take into account a large number of zeros in the coefficient tables
and arithmetic operations on zero components are still performed,
although this is an extra waste of processor time.

One way to achieve versatility and acceptable performance is to
generate code for a numerical method step. This approach minimizes the
number of arithmetic operations and saves memory, since the zero
coefficients of the method do not have to be stored.

We implemented a code generator for the three stochastic numerical
methods mentioned above:\begin{itemize}
\item scalar method with strong convergence $p_s = 1.5$,
\item vector method with strong convergence $p_s = 1.0$,
\item vector method with weak convergence of $p_s = 2.0$.
\end{itemize}
We use Python to implement the code generator and
Jinja2~\cite{l_jinja2} template engine. This template engine was
originally created to generate HTML code, but its syntax is universal
and allows you to generate text of any kind without reference to any
programming or markup language.

Information about the coefficients of each particular method is stored
as a JSON file of the following structure:
\begin{verbatim}
{
  "name": "method's name (the future name of the function)",
  "description": "method's short description",
  "stage": 4,
  "det_order": "2.0",
  "stoch_order": "1.5",
  "A0": [...],
  "B0": [...],
  "A1": [...],
  "B1": [
    ["0", "0", "0", "0"],
    ["1/2", "0", "0", "0"],
    ["-1", "0", "0", "0"],
    ["-5", "3", "1/2", "0"]
  ],
  "c0": ["0", "3/4", "0", "0"],
  "c1": ["0", "1/4", "1", "1/4"],
  "a": ["1/3", "2/3", "0", "0"],
  "b1": ["-1", "4/3", "2/3", "0"],
  "b2": ["-1", "4/3", "-1/3", "0"],
  "b3": ["2", "-4/3", "-2/3", "0"],
  "b4": ["-2", "5/3", "-2/3", "1"]
}
\end{verbatim}

The parameter \texttt{stage} is the number of method's stages,
\texttt{det\_order} is the error order of the deterministic part
($p_d$), \texttt{stoch\_order} is the error order of the stochastic
part ($p_s$), \texttt{name} is the name of the method, which will then
be used to create the name of the generated function, so it should be
written in one word without spaces. All other parameters are the
coefficients of the method. In this case, we give the coefficients of
the scalar method with strong convergence $p_s=1.5$, omitting the
coefficients $\mathbf{a}_{0}$, $\mathbf{a}_{1}$ and $\mathbf{B}_{0}$
to save text space. It is necessary to note that the values of the
coefficients can be specified in the form of rational fractions, for
which they should be presented as JSON strings and enclosed in double
quotes.

For internal representation of stochastic numerical methods we created
three Python classes: \texttt{ScalarMethod},
\texttt{StrongVectorMethod} and \texttt{WeakVectorMethod}. The
implementation of these classes is contained in the file
\texttt{coefficients\_table.py}. The constructors of these classes
read the JSON file and, based on them, create objects, which can later
be used for code generation. The Fraction class from the Python
standard library is used to represent rational coefficients. Each
class has a method that generates a coefficient table in \LaTeX{} format.

The file \texttt{stoch\_rk\_generator.py} is a script which handles
the \texttt{jinja2} templates and, based on them, generates a code of
python functions. For vector stochastic methods, a code is generated
for dimensions up to $6$. Functions are named based on the information
specified in \texttt{JSON} files, such as \texttt{strong\_srk1w2},
\texttt{strong\_srk2w5}, \texttt{weak\_srk2w6}, and so on.

In addition to the code in \texttt{Python}, \LaTeX{} formulas
are generated. It allows one to check the correctness of the
generator. For example, we give below the formula generated
automatically based on the data from \texttt{JSON} file for
Runge--Kutta method \texttt{strong\_srk1w2} with stages $s = 3$, and
$2$ dimensioned Wiener process. Nonzero coefficients of the method are
as follows:
\begin{gather*}
  A_{01}^2 = 1,\; A_{11}^2 = 1,\; A_{11}^3 = 1,\; B_{11}^2 = 1,\; B_{11}^3 = -1,\\
  a_1 = 1/2,\; a_2 = 1/2,\; c_{0}^2 = 1,\; c_{1}^2 = 1,\; c_{1}^3 =
  1,\; b^{1}_1 = 1,\; b^{2}_2 = 1/2,\; b^{2}_3 = -1/2.
\end{gather*}

The numerical scheme formulas are quite cumbersome, despite the large
number of zeros in the coefficient table:
\begin{equation}
  \begin{gathered}
    X^{01\alpha} = x^{\alpha}_{n},\; X^{ 1 1 \alpha} =
    x^{\alpha}_{n},\; X^{ 2 1 \alpha} = x^{\alpha}_{n},
    \\
    X^{02\alpha} = x^{\alpha}_{n} + h_{n}\Big[A_{01}^2 f^{\alpha}(t_n,
    X^{01\beta})\Big],
    \\
    \begin{multlined}
      X^{ 1 2 \alpha} = x^{\alpha}_{n} + h_{n}\Big[A_{11}^2
      f^{\alpha}(t_n, X^{01\beta})\Big]
      \\ {} + %
      B_{11}^2 G^{\alpha}_1(t_n,
      X^{ 1 1 \beta}) \frac{I^{ 1 1 }(h_{n})}{\sqrt{h_{n}}} +
      B_{11}^2 G^{\alpha}_2(t_n, X^{ 2 1 \beta}) \frac{I^{ 2
          1}(h_{n})}{\sqrt{h_{n}}},
    \end{multlined}
    \\
    \begin{multlined}
      X^{ 2 2 \alpha} = x^{\alpha}_{n} + h_{n}\Big[A_{11}^2
      f^{\alpha}(t_n, X^{01\beta})\Big]
      \\ {} + %
      B_{11}^2 G^{\alpha}_1(t_n, X^{
        1 1 \beta}) \frac{I^{ 1 2 }(h_{n})}{\sqrt{h_{n}}}  + B_{11}^2
      G^{\alpha}_2(t_n, X^{ 2 1 \beta}) \frac{I^{ 2 2
        }(h_{n})}{\sqrt{h_{n}}},
    \end{multlined}
    \\
    \begin{multlined}
      X^{ 1 3 \alpha} = x^{\alpha}_{n} + h_{n}\Big[A_{11}^3
      f^{\alpha}(t_n, X^{01\beta})\Big]
      \\ {} + %
      B_{11}^3 G^{\alpha}_1(t_n, X^{
        1 1 \beta}) \frac{I^{ 1 1 }(h_{n})}{\sqrt{h_{n}}}  + B_{11}^3
      G^{\alpha}_2(t_n, X^{ 2 1 \beta}) \frac{I^{ 2 1
        }(h_{n})}{\sqrt{h_{n}}},
    \end{multlined}
    \\
    \begin{multlined}
      X^{2 3 \alpha} = x^{\alpha}_{n} + h_{n}\Big[A_{11}^3
      f^{\alpha}(t_n, X^{01\beta})\Big]
      \\ {} + %
      + B_{11}^3 G^{\alpha}_1(t_n, X^{
        1 1 \beta}) \frac{I^{ 1 2 }(h_{n})}{\sqrt{h_{n}}}
      +  B_{11}^3
      G^{\alpha}_2(t_n, X^{ 2 1 \beta}) \frac{I^{ 2 2
        }(h_{n})}{\sqrt{h_{n}}},
    \end{multlined}
    \\
    \begin{multlined}
      x^{\alpha}_{n+1} = x^{\alpha}_{n} + h_{n} \Big[a_1
      f^{\alpha}(t_{n}, X^{01\beta}) + a_2 f^{\alpha}(t_{n} +
      c_0^2h_{n}, X^{02\beta})\Big] + b^{1}_1
      I^1(h_{n})G^{\alpha}_1(t_{n}, X^{ 1 1 \beta})
      \\ {} + %
      b^{2}_2 \sqrt{h_{n}}G^{\alpha}_1(t_{n} + c_{1}^2 h_{n}, X^{ 1
        2 \beta}) + b^{2}_3 \sqrt{h_{n}}G^{\alpha}_1(t_{n} + c_{1}^3
      h_{n}, X^{ 1 3 \beta})
      \\ {} + %
      b^{1}_1 I^2(h_{n})G^{\alpha}_2(t_{n}, X^{ 2 1 \beta}) + b^{2}_2
      \sqrt{h_{n}}G^{\alpha}_2(t_{n} + c_{1}^2 h_{n}, X^{ 2 2 \beta})
      \\ {} + %
      b^{2}_3 \sqrt{h_{n}}G^{\alpha}_2(t_{n} + c_{1}^3 h_{n}, X^{ 2 3
        \beta}).
    \end{multlined}
  \end{gathered}
\end{equation}

\section{Parallel SDE integration with weak numerical methods}
Stochastic numerical methods with strong convergence are well suited
for computing a specific trajectory of SDE solution. If we are not
interested in a specific trajectory, but in some probabilistic
characteristics (distribution of a random process, mathematical
expectation, variance, etc.), then we should use numerical methods
with weak convergence.

In the case of numerical methods with weak convergence, we have to use
Monte Carlo method. It means that we should solve our SDE system
multiple times and each time with different trajectory. The error of
the Monte Carlo method depends on the number of trials $N$ as
$\sqrt{N}$, so to achieve the accuracy of $10^{-3}$ we need $10^{6}$
trials. However, since the trajectories of the Wiener process are
independent, the SDE for each specific trajectory can be solved
independently in parallel mode.

We have implemented a script in Python, which allows to find solutions
of SDE for $N$ different trajectories in parallel mode by spawning a
given number of processes. For processes spawning we use
\texttt{multiprocessing} module. The following features of the Cpython
interpreter should be noted.
\begin{itemize}
\item Because of the global interpreter lock (GIL), it is not possible
  to use threads for the Monte Carlo method. The standard
  \texttt{threading} module is only suitable for asynchronous tasks.
\item When using processes, you should reinitialize the random number
  generator with new seed for each process separately, because
  otherwise all generated processes will generate the same sequence of
  random numbers.
\end{itemize}

The source code of the implemented script is located in the tests
directory. It is based on two functions.
\begin{itemize}
\item Function \texttt{calculation} performs the necessary
  calculations for a given number of trajectories. As arguments, the
  function takes the drift vector, the diffusion matrix, the required
  number of simulations, the initializing value for the random
  generator, the initial value of the SDU solution, the number of
  steps of the Wiener process, the time interval at which it is
  necessary to carry out integration, the dimension of the Wiener
  process and optionally the function for testing the obtained
  solution for adequacy.
\item Function \texttt{run\_parallel} distributes the Monte Carlo
  tests equally between processes, creates a pool of processes, and
  runs them. Each process performs the function \texttt{calculation}.
\end{itemize}

When carrying out a large number of tests, the storage of all the
resulting trajectories requires a significant amount of
RAM. Therefore, it is more reasonable to immediately decide what
probabilistic characteristics we need and calculate them using on-line
algorithms. For example, to calculate the average trajectory, we use
the following formula
\[
\bar{\mathbf{x}}_{n} = \bar{\mathbf{x}}_{n-1} + \dfrac{\mathbf{x}_n -
  \bar{\mathbf{x}}_{n-1}}{n}.
\]
This formula allows you to update the mean values of all path steps
$\bar{\mathbf{x}}_{n}$ based on the previous mean values
$\bar{\mathbf{x}}_{n-1}$ and the current value $\mathbf{x}_n$. As a
result, each process must store only one array of constant length,
which saves memory.

\section{Conclusion}
\label{sec:conclusion}

Stochastic numerical schemes with convergence order higher than 0.5
are considered. It is shown that such methods are much more
complicated than equivalent numerical methods for systems of ordinary
differential equations. Their specifics makes efficient software
implementation of such methods not a trivial task. We discuss an
approach based on automatic generation of code, which allows to obtain
an efficient implementation of the methods and gives the possibility
to use any table of coefficients. We also give a short description of
our program and a provide url link to the repository with the source
code.

\begin{acknowledgments}

The publication has been prepared with the support of the ``RUDN University Program 5-100'' 
and funded by Russian Foundation for Basic Research (RFBR) according to the research project 
No~16-07-00556.

\end{acknowledgments}

  \bibliography{rk-stoch,self}

\end{document}


\title{Трудности при реализации стохастических численных методов Рунге--Кутты}

\author{М. Н. Геворкян}
\email{gevorkyan-mn@rudn.ru}
\affiliation{Кафедра прикладной информатики и теории вероятностей,\\
  Российский университет дружбы народов,\\
  117198, Москва, ул. Миклухо-Маклая, д. 6}

\author{А. В. Демидова}
\email{demidova-av@rudn.ru}
\affiliation{Кафедра прикладной информатики и теории вероятностей,\\
  Российский университет дружбы народов,\\
  117198, Москва, ул. Миклухо-Маклая, д. 6}

\author{А. В. Королькова}
\email{korolkova-av@rudn.ru}
\affiliation{Кафедра прикладной информатики и теории вероятностей,\\
  Российский университет дружбы народов,\\
  117198, Москва, ул. Миклухо-Маклая, д. 6}

\author{Д. С. Кулябов}
\email{kulyabov-ds@rudn.ru}
\affiliation{Кафедра прикладной информатики и теории вероятностей,\\
  Российский университет дружбы народов,\\
  117198, Москва, ул. Миклухо-Маклая, д. 6}
\affiliation{Лаборатория информационных технологий,\\
  Объединённый институт ядерных исследований,\\
  ул. Жолио-Кюри 6, Дубна, Московская область, Россия, 141980}

\begin{abstract}
  В данной статье рассматриваются стохастические численные методы типа
  Рунге--Кутты с слабой и сильной сходимостями для систем
  стохастических дифференциальных уравнений в форме Ито. Даётся
  краткий обзор основных стохастических численных методов и сведений
  из теории стохастических дифференциальных уравнений. Далее
  излагается и мотивируется подход к реализации данных методов с
  помощью генерации исходного кода. Обсуждаются детали реализации и
  используемые языки программирования и библиотеки.
\end{abstract}

  \keywords{стохастические дифференциальные уравнения, стохастические численные методы,
    автоматическая генерация кода, Python, Julia, шаблонизатор}

\maketitle

\section{Введение}
\label{sec:intro}

При моделировании технических систем с управлением возникает
необходимость исследования их характеристик.
Также необходимо исследование влияния параметров систем на эти
характеристики.
В системах с управлением возникает такое паразитное явление, как
автоколебательный режим.
Нами проводились исследования по определению области
возникновения автоколебаний. Однако параметры этих автоколебаний нами
не исследовались. В данной статье мы предлагаем использовать метод
гармонической линеаризации для данной задачи. Этот метод применяется в
теории управления, однако данный раздел математики достаточно редко
используется в классическом математическом моделировании. Авторы
предлагаю методическую статью, призванную познакомить неспециалистов с
применением этого метода.

\section{Введение}
\label{sec:intro}

В статье авторов~\cite{kulyabov:2016:rk-stochastic} описывалась реализация стохастических
численных методов типа Рунге-Кутты на языке Python~\cite{L_Python} с
использованием библиотек NumPy и SciPy~\cite{L_Scipy}. Выбор языка был
продиктован удобством программирования и возможностью работать с
многомерными массивами как с тензорами (функции \verb|tensor_dot| и
\verb|einsum| из библиотеки NumPy). Однако производительность
созданных функций была на низком уровне, и не столько из-за низкого
быстродействия языка Python, сколько из-за большого числа вложенных
циклов (до семи циклов). В данной статье мы рассматриваем
альтернативных подход к реализации стохастических численных методов,
основанный на генерации кода функций.

Данная статья состоит из трех разделов. В первом разделе дается обзор
основных источников и излагаются сведения из теории стохастических
дифференциальных уравнений (СДУ) и численных методов для их
решения. Во втором разделе приведены стохастические численные схемы
для скалярных СДУ с сильной сходимостью и для систем СДУ с сильной и
слабой сходимостью. Наряду с общими схемами приводятся несколько
таблиц коэффициентов, которые позволяют реализовать уже конкретный
численный метод. Наконец в третьем разделе мотивируется применение
генерации кода для стохастических численных методов и описываются
некоторые детали реализованного нами генератора (на языке Python с
использованием шаблонизатора Jinja2~\cite{L_Jinja2}).

\section{Обзор основных источников}

Первым, кто использовал стохастические броуновские процессы для
математического моделирования был французский математик, ученик
А. Пуанкаре Л. Ж.-Б. А. Башелье (1870--1946) в 1900г. в
работе~\cite{L_Bachelier}.

В данном разделе приведем краткий обзор доступных источников,
посещенных стохастическим методам Рунге-Кутты. Особо отметим, что
такие численные схемы должны быть многостадийными и в них не должно
присутствовать частных производных от $\mathbf{f}(t,\mathbf{X})$ и
$\mathbf{G}(t,\mathbf{X})$.

В качестве учебника по численным методам решения стохастических
дифференциальных уравнений можно указать книгу П. Е. Клодена и
Е. Платена~\cite{L_Kloeden_Platen}. В книге дается краткое введение в
теорию стохастических дифференциальных уравнений Ито и Стратоновича и
области их применений. Последние две трети книги посвящены изложению
численных методов в смысле строгой и слабой аппроксимаций, в том числе
излагается и ряд методов Рунге-Кутты.

В качестве краткого введения в предмет можно использовать обзорную
статью А. В. Лукшина и С. Н. Смирнова~\cite{L_Lukshin_1990}, в которой
тезисно изложены подходы к численному решению стохастических уравнений
по состоянию на 1990 год. В более новой статье~\cite{L_Ereshko_2008}
А. Ф. Ерешко и Д. В. Филатовой также приведены некоторые численные
методы (в том числе и один метод Рунге-Кутты и метод Г. Н.
Милштейна~\cite{L_Milstein_1974, L_Milstein_1979, L_Milstein_1986}), а
также даны результаты численного эксперимента.

Диссертация А. Росслера~\cite{L_Andreas_2003} является
последовательным изложением теории стахостических методов
Рунге-Кутты. Автор рассматривает аппроксимацию систем СДУ Ито и
Стратоновичав в слабом смысле для скалярного и многомерного процесса
Винера. После краткого обзора литературы автор переходит к изложению
теории помеченных деревьев, которая используется для вывода условий
порядка в случае детерминированных методов Рунге-Кутты (смотри,
например,~\cite{L_Butcher_2003,L_Hairer_2008_en,L_Hairer_1990_ru}, и
которая была распространена на случай СДУ. Росслер рассматривает слабо
сходящиеся стохастические методы Рунге-Кутты для систем СДУ Ито и
Стратоновича как для скалярного так и для многомерного процесса
Винера. В третьей и пятой части диссертации приведены конкретные
реализации явных стохастических численных методов Рунге-Кутты (в том
числе и вложенных) для слабой сходимости.

Дальнейшие результаты Росслера изложены в
статьях~\cite{L_Debrabant_2007, L_Debrabant_2013} в соавторстве с
К. Дебрабант. В препринте~\cite{L_Debrabant_2013} продолжена
классификация стохастических методов Рунге-Кутты со слабой сходимостью
на основе условий порядка. Приведено несколько конкретных реализаций и
даны результаты численных экспериментов. В другом
препринте~\cite{L_Rossler_2010} приведены таблицы методов четвертой
стадийности и сильного порядка сходимости $p=3.0$.

Первым стохастическим методом Рунге-Кутты можно считать метод
Эйлера--Маруйамы, изложенный Маруйамой в
статье~\cite{L_Maruyama_1955}. Первым же систематическое исследование
стохастических численных методов Рунге-Кутты сильного порядка
сходимости $p=1.0$ провели В. Рюмелин~\cite{L_Rumelin_1982} и
Е. Платен в диссертации~\cite{L_Platen_1984}.

Большой вклад внесли П. М. Бурраж и К. Бурраж в серии
статей~\cite{L_Burrage_1996,L_Burrage_1997,L_Burrage_1998,L_Burrage_1999,L_Burrage_2000}. В
этих статьях они не только изучили методы сильного порядка $p=1.5$, но
также распространили теорию помеченных деревьев на стохастический
случай.

В статье R. Soheili и M. Namjoo~\cite{L_Soheili_2007} получены три
метода с сильной сходимостью $p=1.0$ и проведено численное сравнение с
методом из книги~\cite{L_Kloeden_Platen}.

Одни из первых методов со слабой сходимостью приведены в
книге~\cite{L_Kloeden_Platen}. Дальнейшее развитие они получили в
статье Й. Комори и Т. Митсуи~\cite{L_Komori_Mitsui_1995} и в статье
\cite{L_Mackevicius_1994}. В статье~\cite{L_Tocino_2001} получены два
трехстадийных метода слабой сходимости $p=2.0$, а также проведены
численные эксперименты.

Обратим внимание на статьи~\cite{L_Wilkie_2004, L_Wilkie_2005} в
которых авторы используют классический метод Рунге--Кутты и разложение
Ито--Тейлора для построения численной схемы для стохастического
дифференциального уравнения, в которую вошли частные производные
вектора сноса $\mathbf{f}(t,\mathbf{X})$ и матрицы диффузии
$\mathbf{G}(t,\mathbf{X})$. Такие схемы мы в данной статье
рассматривать не будем.

В виду крайней сложности дальнейшего повышения порядка точности
стохастических численных схем, современные работы в основном
концентрируются на получении численных схем для частных случаев
СДУ. Можно выделить работы посвященные симплектическим стохастическим
численным методам
Рунге--Кутта~\cite{L_Burrage2012,L_Ma2015,L_Zhou2017} и неявным
схемам, в частности стохастическим аналогам метода
Розенброка~\cite{L_Amiri2017}.

Из литературы на русском языке стоит выделить
книгу~\cite{L_Kuznezov2010} в которой дан обзор большого количества
математических моделей на основе стохастических систем уравнений из
различных областей естественных наук. Также автор рассматривает вопрос
аппроксимации интеграла Ито, и построение ряда численных схем с
сильной сходимостью на основе этих аппроксимаций.

\section{Случайный процесс Винера  и программное генерирование его траекторий}
Случайный процесс $W(t)$, $t \geqslant 0$ называется скалярным
\emph{процессом Винера} если выполняются следующие
условия~\cite{L_Oksendal_ru,L_Kloeden_Platen}:
  \begin{itemize}
    \item $\mathrm{P}\{W(0)=0\} = 1$, иначе говоря $W(0) = 0$ почти наверное;
    \item $W(t)$ --- процесс с независимыми приращениями, то есть $\{\Delta W_{i}\}^{N-1}_{0}$ --- независимые случайные величины: $\Delta W_{i} = W(t_{i+1}) - W(t_{i})$ и $0 \leqslant t_{0} < t_{1} < t_{2} < \ldots  < t_{N} \leqslant T$;
    \item $\Delta W_{i} = W(t_{i+1}) - W(t_{i}) \sim \mathcal{N}(0,t_{i+1}-t_{i})$ где $0\leqslant t_{i+1} < t_{i} < T$, $i=0,1,\ldots,N-1$
  \end{itemize}

  Обозначение $\Delta W_{i} \sim \mathcal{N}(0,\Delta t_{i})$ говорит
  о том, что $\Delta W_{i}$ --- нормально распределенная случайная
  величина с математическим ожиданием
  $\mathbb{E}[\Delta W_{i}] = \mu = 0$ и дисперсией
  $\mathbb{D}[\Delta W_{i}] = \sigma^{2} = \Delta t_{i}$.

  Винеровский процесс является моделью \emph{броуновского движения}
  (хаотического блуждания). Если рассмотреть процесс $W(t)$ в те
  моменты времени
  $0 = t_{0} < t_{1} < t_{2} < \ldots < t_{N-1} < t_{N}$ когда он
  испытывает случайные аддитивные изменения, то непосредственно из
  определения следует:
\[
    W(t_{1}) = W(t_{0}) + \Delta W_0, W(t_{2}) = W(t_{1}) + \Delta W_1, \ldots, W(t_{N}) = W(t_{N-1}) + \Delta W_{N-1},
\]
где $\Delta W_{i} \sim \mathcal{N}(0,\Delta t_{i})$, $\forall i = 0,\ldots,N-1$.

Если величину $W(t_{N})$ записать в виде накопленной суммы приращений:
\[
  W(t_{n}) = W(t_{0}) + \sum\limits_{i=0}^{n}\Delta W_{i},\;\;n=0,\ldots,N-1,
\]
и учесть, что $\mathbb{E}[\Delta W_{i}] = 0$ и $\mathbb{D}[\Delta W_{i}] = \Delta t_{i}$, то можно показать, что сумма нормально распределенных случайных чисел $\Delta W_{i}$ также является нормально распределенным случайным числом:
\[
  \mathbb{E}\sum\limits_{i=0}^{N}\Delta W_{i} = 0, \;\; \mathbb{D}\sum\limits_{i=0}^{n}\Delta W_{i} = \sum\limits_{i=0}^{n} \Delta t_{i} = t_{n} - t_{0},\;\; \sum\limits_{i=0}^{n}\Delta W_{i} \sim \mathcal{N}(0,t_{n}-t_{0}).
\]

Многомерный винеровский процесс
$\mathbf{W}(t)\colon \Omega\times[t_{0},T]\to \mathbb{R}^{m}$
определяют как случайный процесс составленный из совместно независимых
одномерных винеровских процессов
$W^{1}(t),\ldots,W^{m}(t)$. Приращения
$\Delta W^{\alpha}_{i},\;\forall \alpha = 1,\ldots,m$ являются
совместно независимыми нормально распределенными случайными
величинами. С другой стороны вектор $\Delta W^{\alpha}_{i}$ можно
представить как многомерную нормально распределенную случайную
величину с вектором математических ожиданий $\mu = \mathbf{0}$ и
диагональной матрицей ковариаций.

В случае многомерного случайного процесса следует сгенерировать уже $m$ последовательностей из $N$ нормально распределенных случайных величин.

\section{Стохастические интегралы и СДУ для скалярного процесса Винера}

Формулировка стохастического дифференциального уравнения начинается с
введения операции стохастического интегрирования.

Пусть $g(t,x(t))$ непрерывная скалярная функция, зависящая от
случайного процесса $x(t)$, $W(t)$ --- процесс Винера,
$t\in[t_{0},T]$,
$0 \leqslant t_{0} < t_{1} < t_{2} < \ldots < t_{n-1} < t_{n}
\leqslant T < \infty$ --- разбиение отрезка $[t_{0},T]$, тогда
  \[
    I(\theta) = \int\limits_{t_{0}}^{T}g(t,x(t))\mathrm{d} W(t) = \mathop{\mathrm{l.i.m.}\;}\limits_{n\to\infty}\sum^{n-1}_{i=0}g(\theta t_{i+1} + (1-\theta)t_{i}, \theta x(t_{i+1}) + (1-\theta)x(t_{i}))(W(t_{i+1}) - W(t_{i}))
  \]
где $\theta\in[0,1]$. $I(\theta)$ называется \emph{стохастическим
  интегралом}. Строгое обоснование существования данного интеграла для
широкого класса функций см.~\cite[глава
3]{L_Oksendal_ru}.

В в физике и прикладной математике нашли применение два частных случая
стохастического интеграла:
\begin{itemize}
  \item \emph{интеграл Ито} (назван в честь японского математика К. Ито)
  \[
    I = I(0) = \mathop{\mathrm{l.i.m.}\;}\limits_{n\to\infty}\sum^{n-1}_{i=0}g(t_{i}, x(t_{i}))(W(t_{i+1}) - W(t_{i})),
  \]
  \item \emph{интеграл Стратоновича} (назван в честь советского физика Р. Л. Стратоновича)
  \[  
    I_{0.5} = I(0.5) = \mathop{\mathrm{l.i.m.}\;}\limits_{n\to\infty}\sum^{n-1}_{i=0}g\left( \frac{t_{i+1}+t_{i}}{2},\frac{x(t_{i+1}) + x(t_{i})}{2}\right)(W(t_{i+1}) - W(t_{i})).
  \]
\end{itemize}

После введения операции стохастического интегрирования, можно записать интегрального уравнения для случайного процесса $x(t)$~\cite{L_Oksendal_ru}.
\[
  x(t_{k}) = x(t_{0}) + \int\limits_{0}^{t}f(\tau,x(\tau))\mathrm{d} \tau + \underbrace{\int\limits_{0}^{t}g(\tau,x(\tau))\mathrm{d} W.}_{\text{стохастический интеграл}}
\]
Вышеприведенному интегральному уравнению соответствует стохастическое дифференциальное уравнение (СДУ):
\[
  \mathrm{d} x(t) = f(t,x(t))\mathrm{d} t + g(t, x(t))\mathrm{d} W.
\]
При этом необходимо помнить, что дифференциал типа $\mathrm{d} x$ не
является обычным <<малым>> приращением функции $x(t)$, а является
случайной величиной. Фактически СДУ символической записью предела
итерационной схемы. Процесс Винера $W(t)$, приращения которого входят
в состав СДУ, называется ведущим процессом для данного СДУ (driving
process).

\section{СДУ в форме Ито для многомерного процесса Винера}

Сформулируем определение системы стохастических дифференциальных
уравнений. Введем вероятностное пространство
$(\Omega, \mathscr{A}, \mathbb{P})$, где $\Omega$ --- пространство
элементарных событий, $\mathscr{A}$ --- сигма-алгебра подмножеств
пространства $\Omega$, $\mathbb{P}$ --- вероятностная мера. На отрезке
$[t_{0}, T] \in \mathbb{R}^{1}$ определена переменная $t$ имеющая
физический смысл времени.

Рассмотрим случайный процесс
$\mathbf{x}(t) = (x^{1}(t), \ldots, x^{d}(t))^{T}$, где
$\mathbf{x}(t)$ принадлежит функциональному пространству
$\mathrm{L}^{2}(\Omega)$ с нормой $\| \cdot \|$. Будем считать, что
случайный процесс $\mathbf{x}(t)$ является решением для
стохастического дифференциального уравнения (СДУ) в форме
\emph{Ито}~\cite{L_Oksendal_ru, L_Kloeden_Platen}:
\[
  \mathbf{x}(t) = \mathbf{f}(t,\mathbf{x}(t))\mathrm{d} t + \mathbf{G}(t,\mathbf{x}(t))\mathrm{d} \mathbf{W},
\]
где $\mathbf{W} = (W^{1},\dots,W^{m})^T$ --- многомерный винеровский
процесс, называемый \emph{ведущим (driver)} процессом СДУ.  Функция
$\mathbf{f}\colon[t_{0}, T]\times\mathbb{R}^{d} \to \mathbb{R}^{d}$
называется \emph{вектором сноса}, а матричнозначная функция
$\mathbf{G}\colon[t_{0}, T]\times\mathbb{R}^{d}\times\mathbb{R}^{m}
\to \mathbb{R}^{d}\times\mathbb{R}^{m}$ --- \emph{матрицей
  диффузии}. Кроме того
$\mathbf{f}(t,\mathbf{x}(t)) =
(f^{1}(t,\mathbf{x}),\ldots,f^{d}(t,\mathbf{x}))^{T}$, а матрица
$\mathbf{G}$ имеет вид:
\[
  \mathbf{G} = 
  \begin{bmatrix}
    g^{1}_{1}(t,\mathbf{x}) & g^{1}_{2}(t,\mathbf{x}) & \ldots & g^{1}_{m}(t,\mathbf{x})\\
    g^{2}_{1}(t,\mathbf{x}) & g^{2}_{2}(t,\mathbf{x}) & \ldots & g^{2}_{m}(t,\mathbf{x})\\
    \vdots & \vdots & \ddots & \vdots \\
    g^{d}_{1}(t,\mathbf{x}) & g^{d}_{2}(t,\mathbf{x}) & \ldots & g^{d}_{m}(t,\mathbf{x})
  \end{bmatrix}
\]
То же уравнение можно переписать в индексном виде:
\[
  x^{\alpha}(t) = f^{\alpha}(t,x^{\gamma}(t))\mathrm{d} t + \sum\limits_{\beta = 1}^{m}g^{\alpha}_{\beta}(t,x^{\gamma}(t))\mathrm{d} W^{\beta},
\]
где $\alpha,\gamma = 1,\ldots,d$, $\beta = 1,\ldots,m$, а $f^{\alpha}(t,x^{\gamma}(t)) = f^{\alpha}(t,x^{1}(t),\ldots,x^{d}(t))$.

На отрезке $[t_{0}, T]$ введем сетку $t_{0}<t_{1}<\ldots<t_{N}=T$ с
шагом $h_{n} = t_{n+1} - t_{n}$, где $n = 0,\ldots,N-1$ и максимальный
шаг сетки $h = \max{\{h_{n-1}\}^{N}_{1}}$. Далее будем полагать, что
сетка равномерная, то $h_{n} = h = \mathrm{const}$. $\mathbf{x}_{n}$
--- сеточная функция, аппроксимирующая случайный процесс
$\mathbf{x}(t)$, так что $\mathbf{x}_{0} = \mathbf{x}(t_{0})$,
$\mathbf{x}_{n} \approx \mathbf{x}(t_{n})\;\forall n = 1,\ldots,N$.

\section{Вычисление и аппроксимация кратных интегралов Ито специального вида}

Здесь мы не будем вдаваться в общую теорию кратных стохастических
интегралов Ито, за которой отсылаем читателя к
книге~\cite{L_Kloeden_Platen}. Здесь же рассмотрим кратные интегралы
частного вида, которые входят в формулы стохастических численных
методов.

В общем случае для конструирования численных схем, порядок сходимости
которых был бы больше чем $p=\frac12$, необходимо включение в формулы
этих схем однократных, двукратных и трехкратных интегралов Ито
следующего вида:
\[
  I^{\alpha}(t_{n}, t_{n+1}) = I^{\alpha}(h_{n}) = \int\limits_{t_{n}}^{t_{n+1}} \mathrm{d}W^{\alpha}(\tau),
\]
\[
  I^{\alpha\beta}(t_{n}, t_{n+1}) = I^{\alpha\beta}(h_{n}) = \int\limits_{t_{n}}^{t_{n+1}}\int\limits_{t_{n}}^{\tau_{1}}\mathrm{d}W^{\alpha}(\tau_{2})\mathrm{d}W^{\beta}(\tau_{1}),
\]
\[
  I^{\alpha\beta\gamma}(t_{n}, t_{n+1}) = I^{\alpha\beta\gamma}(h_{n}) = \int\limits_{t_{n}}^{t_{n+1}}\int\limits_{t_{n}}^{\tau_{1}}\int\limits_{t_{n}}^{\tau_{2}}\mathrm{d}W^{\alpha}(\tau_{3})\mathrm{d}W^{\beta}(\tau_{2})\mathrm{d}W^{\gamma}(\tau_{1}),
\]
где $\alpha,\beta,\gamma=0\ldots,m$ и $W^{\alpha},\alpha=1,\ldots,m$
--- компоненты многомерного винеровского процесса. В случае
$\alpha,\beta,\gamma = 0$, приращение $\mathrm{d} W^{0}(\tau)$
полагается равным $\mathrm{d} \tau$.

Возникает задача выражения этих интегралов через приращения
$\Delta W^{i}_{n} = W^{i}(t_{n+1}) - W^{i}(t_{n})$. Несмотря на
кажущуюся простоту, сделать это можно далеко не для всех возможных
комбинаций индексов. Рассмотрим в начале те случаи, когда
аналитическое выражение получить возможно, а затем перейдем к тем
случаям, когда приходится прибегать к аппроксимирующим формулам.

В случае однократного интеграла задача тривиальна и аналитическое
выражение можно получить для любого индекса $\alpha$:
\[
  I^{0}(h_{n}) = \Delta t_{n} = h_{n},\;\; I^{\alpha}(h_{n}) = \Delta W^{\alpha}_{n},\;\alpha=1,\ldots,m.
\]

В случае двукратного интеграла $I^{\alpha\beta}(h_{n})$ точная формула
имеет место лишь при $\alpha=\beta$:
\[
  I^{00}(h_{n}) = \dfrac{1}{2}\Delta t_{n} = \dfrac{1}{2}h^{2}_{n},\;\; I^{\alpha\alpha}(h_{n}) = \dfrac{1}{2}\left( (\Delta W^{\alpha}_{n})^2 - \Delta t_{n} \right),\;\alpha=1,\ldots,m,
\]
в остальных же случаях при $\alpha\neq \beta$ выразить
$I^{\alpha\beta}(h_{n})$ через приращения $\Delta W^{\alpha}_{n}$ и
$\Delta t_{n}$ в конечном виде не представляется возможным, поэтому
остается лишь использовать численную аппроксимацию.

Для смешанного случая $I^{0\alpha}$ и $I^{\alpha0}$
в~\cite{L_Andreas_2003} приведены простые формулы следующего вида:
\[
  I^{0\alpha}(h_{n}) = \dfrac{1}{2}h_{n}\left(I^{\alpha}(h_{n}) - \frac{1}{\sqrt{3}}\zeta^{\alpha}(h_{n})\right),
\]
\[
  I^{\alpha0}(h_{n}) = \dfrac{1}{2}h_{n}\left(I^{\alpha}(h_{n}) + \frac{1}{\sqrt{3}}\zeta^{\alpha}(h_{n})\right),
\]
где $\zeta^{\alpha}_{n}\sim \mathcal{N}(0,h_{n})$ многомерные
нормально распределенные случайные величины.

Для общего случая $\alpha,\beta = 1,\ldots,m$ в
книге~\cite{L_Kloeden_Platen} приведены следующие формулы для
аппроксимации двукратного интеграла Ито $I^{\alpha\beta}$:
\begin{gather*}
  I^{\alpha\beta}(h_{n}) = \dfrac{\Delta W^{\alpha}_{n} \Delta W^{\beta}_{n} - h_{n}\delta^{\alpha\beta}}{2} + A^{\alpha\beta}(h_{n}),\\
  A^{\alpha\beta}(h_{n}) = \dfrac{h}{2\pi}\sum\limits^{\infty}_{k=1}\dfrac{1}{k}\left[V^{\alpha}_{k}\left(U^{\beta}_{k} + \sqrt{\frac{2}{h_{n}}} \Delta W^{\beta}_{n}\right) - V^{\beta}_{k}\left(U^{\alpha}_{k} + \sqrt{\frac{2}{h_{n}}}\Delta W^{\alpha}_{n}\right)\right],
\end{gather*}
$V^{\alpha}_{k} \sim \mathcal{N}(0,1)$,
$U^{\alpha}_{k} \sim \mathcal{N}(0,1)$,
$\alpha=1,\ldots,m;\; k=1,\ldots,\infty$; $n=1,\ldots,N$ --- номер
шага численной схемы. Из формул видно, что в случае $\alpha=\beta$
получаем конечное выражение для $I^{\alpha\beta}$, приведенное нами
выше. В случае $\alpha\neq \beta$ приходится суммировать бесконечный
ряд $A^{\alpha\beta}$. Данный алгоритм дает ошибку аппроксимации
порядка $O(h^2/n)$, где $n$ --- число оставленных слагаемых
бесконечного ряда $A^{ij}$.

В статье~\cite{L_Wiktorsson_2001} предложена матричный вид
аппроксимирующих формул. Пусть $\mathbf{1}_{m\times m}$,
$\mathbf{0}_{m\times m}$ --- единичная и нулевая матрицы $m\times m$,
тогда
\begin{gather*}
  \mathbf{I}(h_{n}) = \dfrac{\Delta \mathbf{W}_{n} \Delta \mathbf{W}^{T}_{n} - h_{n}\mathbf{1}_{m\times m}}{2} + \mathbf{A}(h_{n}),\\
  \mathbf{A}(h_{n}) = \dfrac{h}{2\pi}\sum\limits^{\infty}_{k=1}\dfrac{1}{k}\left(\mathbf{V}_{k}(\mathbf{U}_{k} + \sqrt{2/h_{n}} \Delta \mathbf{W}_{n})^{T} - (\mathbf{U}_{k} + \sqrt{2/h_{n}}\Delta \mathbf{W}_{n})\mathbf{V}^{T}_{k}\right),
\end{gather*}
где $\Delta \mathbf{W}_{n}, \mathbf{V}_{k}, \mathbf{U}_{k}$ ---
независимые нормально распределенные многомерные случайные величины:
\begin{gather*}
\Delta \mathbf{W}_{n} = (\Delta W^{1}_{n}, \Delta W^{2}_{n}, \ldots, \Delta W^{m}_{n})^{T} \sim \mathcal{N}(\mathbf{0}_{m\times m}, h_{n}\mathbf{1}_{m\times m}),\\
\mathbf{V}_{k} = (V^{1}_{k}, V^{2}_{k},\ldots, V^{m}_{k})^{T} \sim \mathcal{N}(\mathbf{0}_{m\times m}, \mathbf{1}_{m\times m}),\;\;
\mathbf{U}_{k} = (U^{1}_{k}, U^{2}_{k},\ldots, U^{m}_{k})^{T} \sim \mathcal{N}(\mathbf{0}_{m\times m}, \mathbf{1}_{m\times m}).
\end{gather*}
В случае, если язык программирования поддерживает векторизацию
действий с многомерными массивами, данные формулы могут обеспечить
выигрыш к производительности программы.

Наконец рассмотрим трехкратный интеграл. В единственной численной
схеме, в которой он встречается, необходимо уметь вычислять только
случай одинаковых индексов $\alpha=\beta=\gamma$. Для данного случая
в~\cite{L_Andreas_2003} приведена следующая формула:
\[
  I^{\alpha\alpha\alpha}(h_{n}) = \dfrac{1}{6}\left((I^{\alpha}(h_{n}))^{3} - 3 I^{0}(h_{n})I^{\alpha}(h_{n})\right) = \dfrac{1}{6}\left((\Delta W^{\alpha}_{n})^{3} - 3 h_{n}\Delta W^{\alpha}_{n}\right).
\]
\section{Слабая и сильная сходимости аппроксимирующей функции}
Прежде чем переходить к формулировке численных схем, необходимо
определить критерий точности аппроксимации моделируемого процесса
$\mathbf{x}(t)$ аппроксимирующей сеточной функцией
$\mathbf{x}_{n}$. Таких критериев принято выделять два: \emph{слабый}
и \emph{сильный}.

Последовательность сеточных аппроксимирующих функций
$\{\mathbf{x}_{n}\}^{N}_{1}$ сходится в \emph{сильном смысле} с
порядком $p$ к решению $\mathbf{x}(t)$ СДУ в момент времени $T$ если
существует константа $C>0$ и $\delta_{0}>0$ такая что
$\forall h \in (0,\delta_{0}]$ выполняется условие
\[
\mathbb{E}(\|\mathbf{x}(T) - \mathbf{x}_{N}\|) \leqslant Ch^{p}.
\]

Последовательность сеточных аппроксимирующих функций
$\{\mathbf{x}_{n}\}^{N}_{1}$ сходится в \emph{слабом смысле} с
порядком $p$ к решению $\mathbf{x}(t)$ СДУ в момент времени $T$ если
существует константа $C_{F}>0$ и $\delta_{0}>0$ такая что
$\forall h \in (0,\delta_{0}]$ выполняется условие
\[
\left|\mathbb{E}\left[F(\mathbf{x}(T))\right] -
  \mathbb{E}\left[F(\mathbf{x}_{N})\right]\right| \leqslant
C_{F}h^{p}.
\]

Здесь $F \in C^{2(p+1)}_{\mathrm{P}}(\mathbb{R},\mathbb{R}^{d})$ ---
непрерывный дифференцируемый функционал с полиномиальным ростом.

Если матрица $\mathbf{G}$ обращается в ноль, то условие сильной
сходимости равносильно условию сходимости для детерминированного
случая. В отличие от детерминированного случая порядок сильной
сходимости не обязательно является натуральным числом и может
принимать дробно-рациональные значения.

Важно отметить, что выбор типа сходимости зависит от решаемой
задачи. Увеличение порядка строгой сходимости ведет к увеличению
точности аппроксимации траекторий решения $\mathbf{x}(t)$ СДУ. Если же
требуется вычислить, например, момент случайного процесса
$\mathbf{x}(t)$ или обобщенный функционал вида
$\mathbb{E}[F(\mathbf{x}(t))]$, то следует повышать порядок слабой
сходимости.

\section{Стохастические методы Рунге--Кутты}

\subsection{Численный метод Эйлера--Маруямы}

Простейшим численным методом для решения как скалярных уравнений, так
и систем СДУ, является метод Эйлера--Маруямы, названный так в честь
Гиширо Маруямы (Gisiro Maruyama), который распространил классический
метод Эйлера для ОДУ на случай СДУ~\cite{L_Maruyama_1955}. Метод легко
обобщается на случай многомерного винеровского процесса.
\begin{align*}
  &x^{\alpha}_{0} = x^{\alpha}(t_{0}),\\
  &x^{\alpha}_{n+1} = x^{\alpha}_{n} + f^{\alpha}(t_{n}, x^{\alpha}_{n})h_{n} + \sum\limits^{d}_{\gamma=1}G^{\alpha}_{\beta}(t_{n},x^{\gamma}_{n})\Delta W^{\beta}_{n}.
\end{align*}
Как видно из формул на каждом шаге для вычисления следующего значения
аппроксимирующей функции требуется лишь соответствующее данному шагу
приращение процесса $\Delta W^{\beta}_{n}$. Метод имеет сильный
порядок $(p_{d},p_{s}) = (1.0, 0.5)$. Величина $p_{d}$ обозначает
порядок точности детерминированной части численного метода, то есть
той точности, которую будет давать численный метод при применении его
к уравнению функцией $G(t,x^{\alpha}(t))\equiv 0$. Величина $p_{s}$
обозначает порядок приближения стохастической части уравнения.

\subsection{Стохастические методы Рунге--Кутты сильной сходимости $p=1.5$ для одномерного винеровского процесса}

В случае скалярного СДУ вектор сноса $f^{\alpha}(t, x^{\gamma})$ и
матрица диффузии $G^{\alpha}_{\beta}(t, x^{\gamma})$ вырождаются в
скалярные функции $f(t, x)$ и $g(t, x)$, а ведущий винеровский процесс
$W^{\beta}_t$ становится одномерным $W_t$. Для такого уравнения можно
составить численную схему со сильной сходимостью $p=1.5$:
\begin{align*}
  &X^{i}_{0} = x_{n} + \sum\limits^{s}_{j=1}A_{0j}^{i}f(t_{n} + c^{j}_{0}h_{n}, X^{j}_{0})h_{n} + \sum\limits^{s}_{j=1}B^{i}_{0j}g(t_{n} + c^{j}_{1}h_{n}, X^{j}_{1})\dfrac{I^{10}(h_{n})}{\sqrt{h_{n}}},\\
  &X^{i}_{1} = x_{n} + \sum\limits^{s}_{j=1}A_{1j}^{i}f(t_{n} + c^{j}_{0}h_{n}, X^{j}_{0})h_{n} + \sum\limits^{s}_{j=1}B^{i}_{1j}g(t_{n} + c^{j}_{1}h_{n}, X^{j}_{1})\sqrt{h_{n}},\\
  &x_{n+1} = x_{n} + \sum\limits^{s}_{i=1}a_{i}f(t_{n}+c^{i}_{0}h_{n}, X^{i}_{0})h_{n} + \\
  & \qquad + \sum\limits^{s}_{i=1}\left(b^{1}_{i}I^{1}(h_{n}) + b^{2}_{i}\dfrac{I^{11}(h_{n})}{\sqrt{h_{n}}} + b^{3}_{i}\dfrac{I^{10}(h_{n})}{h_{n}} + b^{4}_{i}\dfrac{I^{111}(h_{n})}{h_{n}}\right)g(t_{n} + c^{i}_{1}h_{n},X^{i}_{1}),
\end{align*}
где $i,j=1,\ldots,s$ ($s$ --- стадийность численного
метода). Обобщенная таблица Бутчера данной схемы имеет
вид~\cite{L_Rossler_2010}:
\[
  {\renewcommand{\arraystretch}{1.5}%
  \begin{array}{c|c|c|c}
    c^{i}_{0} & A^{i}_{0j} & B^{i}_{0j} &\\ \hline
    c^{i}_{1} & A^{i}_{1j} & B^{i}_{1j} &\\ \hline
    & a_{i} & b^{1}_{i} & b^{2}_{i}\\ \hline
    &  & b^{3}_{i} & b^{4}_{i}\\
  \end{array}}
\]

Заметим, что винеровский процесс и его приращения присутствуют в
вышеприведенных формулах в неявном виде кратных интегралов Ито:
$I^{10}(h_{n})$, $I^{1}(h_{n})$, $I^{11}(h_{n})$, $I^{111}(h_{n})$,
общие формулы для которых приведены выше. Для скалярного случая они
упрощаются:
\begin{align*}
  & I^{1}(h_{n}) = \Delta W_{n},\\
  & I^{10}(h_{n}) = \dfrac{1}{2}h_{n}(\Delta W_{n} + \zeta_{n}/\sqrt{3}),\\
  & I^{11}(h_{n})  = \dfrac{1}{2}((\Delta W_{n})^2 - h_{n}), \\
  & I^{111}(h_{n}) = \dfrac{1}{6}((\Delta W_{n})^3 - h_{n}\Delta W_{n}),
\end{align*}
где $\zeta_{n}\sim \mathcal{N}(0,h_{n}), \sigma = \sqrt{h_{n}}$.  В
препринте Росслера~\cite{L_Rossler_2010} для данной схемы приведены
две таблицы Бутчера для метода четвертой стадийности $s=4$.
\[
  \text{\texttt{SRK1W1}: }
  {\renewcommand{\arraystretch}{1.1}%
\begin{array}{c|cccc|cccc|cccc}
  0 & 0 & 0 & 0 & 0 & 0 & 0 & 0 & 0 &  &  &  & \\
  3/4 & 3/4 & 0 & 0 & 0 & 3/2 & 0 & 0 & 0 &  &  &  & \\
  0 & 0 & 0 & 0 & 0 & 0 & 0 & 0 & 0 &  &  &  & \\
  0 & 0 & 0 & 0 & 0 & 0 & 0 & 0 & 0 &  &  &  & \\
  \hline
  0 & 0 & 0 & 0 & 0 & 0 & 0 & 0 & 0 &  &  &  & \\
  1/4 & 1/4 & 0 & 0 & 0 & 1/2 & 0 & 0 & 0 &  &  &  & \\
  1 & 1 & 0 & 0 & 0 & -1 & 0 & 0 & 0 &  &  &  & \\
  1/4 & 0 & 0 & 1/4 & 0 & -5 & 3 & 1/2 & 0 &  &  &  & \\
  \hline
  & 1/3 & 2/3 & 0 & 0 & -1 & 4/3 & 2/3 & 0 & -1 & 4/3 & -1/3 & 0\\
  \hline
  &   &   &   &   & 2 & -4/3 & -2/3 & 0 & -2 & 5/3 & -2/3 & 1
\end{array}}
\]
\[
  \text{\texttt{SRK2W1}: }
  {\renewcommand{\arraystretch}{1.1}%
  \begin{array}{c|cccc|cccc|cccc}
  0 & 0 & 0 & 0 & 0 & 0 & 0 & 0 & 0 &  &  &  & \\
  1 & 1 & 0 & 0 & 0 & 0 & 0 & 0 & 0 &  &  &  & \\
  1/2 & 1/4 & 1/4 & 0 & 0 & 1 & 1/2 & 0 & 0 &  &  &  & \\
  0 & 0 & 0 & 0 & 0 & 0 & 0 & 0 & 0 &  &  &  & \\
  \hline
  0 & 0 & 0 & 0 & 0 & 0 & 0 & 0 & 0 &  &  &  & \\
  1/4 & 1/4 & 0 & 0 & 0 & -1/2 & 0 & 0 & 0 &  &  &  & \\
  1 & 1 & 0 & 0 & 0 & 1 & 0 & 0 & 0 &  &  &  & \\
  1/4 & 0 & 0 & 1/4 & 0 & 2 & -1 & 1/2 & 0 &  &  &  & \\
  \hline
  & 1/6 & 1/6 & 2/3 & 0 & -1 & 4/3 & 2/3 & 0 & -1 & -4/3 & 1/3 & 0\\
  \hline
  &   &   &   &   & 2 & -4/3 & -2/3 & 0 & -2 & 5/3 & -2/3 & 1
\end{array}}
\]
Численную схему реализуемую первой таблицей обозначим как
\texttt{SRK1W1}, а вторую как \texttt{SRK2W2}. Метод \texttt{SRK1W1}
имеет сильный порядок $(p_{d},p_{s}) = (2.0, 1.5)$, а метод
\texttt{SRK2W1} сильный порядок $(p_{d},p_{s}) = (3.0, 1.5)$. Еще один
метод с сильным порядком $p_{s}=1.0$ приведен в
книге~\cite{L_Kloeden_Platen} и его таблица Бутчера имеет вид:
\[
  \text{\texttt{KlPl}: }
  {\renewcommand{\arraystretch}{1.1}%
  \begin{array}{c|cc|cc|cc}
  0&0&0&0&0&&\\
  0&0&0&0&0&&\\ \hline
  0&0&0&0&0&&\\
  0&1&0&1&0&&\\ \hline
  0&1&0&1&0&-1&1\\\hline
   & & &1&0&0&0\\
  \end{array}}
\]


\subsection{Стохастические методы Рунге--Кутты сильной сходимости $p=1.0$ для многомерного винеровского потока}
Для системы СДУ Ито с многомерным винеровским процессом можно
построить стохастическую численную схему Рунге-Кутты сильного порядка
$p_{s} = 1.0$ с использованием однократных и двукратных интегралов
Ито~\cite{L_Rossler_2010}.
\begin{align*}
  &X^{0i\alpha} = x^{\alpha}_{n} + \sum\limits^{s}_{j=1}A^{i}_{0j}f^{\alpha}(t_{n} + c^{j}_{0}h_{n}, X^{0j\beta})h_{n} + \sum\limits_{l=1}^{m}\sum\limits^{s}_{j=1}B^{i}_{0j}G^{\alpha}_{l}(t_{n} + c^{j}_{1}h_{n},X^{lj\beta})I^{l}(h_{n}),\\
  &X^{ki\alpha} = x^{\alpha}_{n} + \sum\limits^{s}_{j=1}A^{i}_{1j}f^{\alpha}(t_{n} + c^{j}_{0}h_{n},X^{0j\beta})h_{n} + \sum\limits_{l=1}^{m}\sum\limits^{s}_{j=1}B^{i}_{1j}G^{\alpha}_{l}(t_{n} + c^{j}_{1}h_{n},X^{lj\beta})\dfrac{I^{lk}(h_{n})}{\sqrt{h_{n}}},\\
  &x^{\alpha}_{n+1} = x^{\alpha}_{n} + \sum\limits^{s}_{i=1}a_{i}f^{\alpha}(t_{n} + c^{i}_{0}h_{n},X^{0i\beta})h_{n} + \sum\limits^{m}_{k=1}\sum\limits^{s}_{i=1}(b^1_{i}I^{k}(h_{n}) + b^2_{i}\sqrt{h_{n}})G^{\alpha}_{k}(t_{n} + c^{i}_{1}h_{n},X^{ki\beta}),
\end{align*}
$n=0,1,\ldots,N-1$; $i=1,\ldots,s$; $\beta,k=1,\ldots,m$; $\alpha = 1,\ldots,d$. Обобщенная таблица Бутчера имеет вид~\cite{L_Rossler_2010}:
\[
  {\renewcommand{\arraystretch}{1.4}%
  \begin{array}{c|c|c|c}
    c^{i}_0 & A^{i}_{0j} & B^{i}_{0j} &\\ \hline
    c^{i}_1 & A^{i}_{1j} & B^{i}_{1j} &\\ \hline
    & a_{i} & b^{1}_{i} & b^{2}_{i} \\
  \end{array}}
\]
В препринте Росслера~\cite{L_Rossler_2010} указаны две таблицы Бутчера
для метода третей стадийности $s=3$:
\[
  \text{\texttt{SRK1Wm}: }
  {\renewcommand{\arraystretch}{1.1}%
  \begin{array}{c|ccc|ccc|ccc}
    0 & 0 & 0 & 0 & 0 & 0 & 0 &  &  & \\
    0 & 0 & 0 & 0 & 0 & 0 & 0 &  &  & \\
    0 & 0 & 0 & 0 & 0 & 0 & 0 &  &  & \\
    \hline
    0 & 0 & 0 & 0 & 0 & 0 & 0 &  &  & \\
    0 & 0 & 0 & 0 & 1 & 0 & 0 &  &  & \\
    0 & 0 & 0 & 0 & -1 & 0 & 0 &  &  & \\
    \hline
    & 1 & 0 & 0 & 1 & 0 & 0 & 0 & 1/2 & -1/2\\
  \end{array}}
\]
\[
  \text{\texttt{SRK2Wm}: }
  {\renewcommand{\arraystretch}{1.1}%
    \begin{array}{c|ccc|ccc|ccc}
    0 & 0 & 0 & 0 & 0 & 0 & 0 &  &  & \\
    1 & 1 & 0 & 0 & 0 & 0 & 0 &  &  & \\
    0 & 0 & 0 & 0 & 0 & 0 & 0 &  &  & \\
    \hline
    0 & 0 & 0 & 0 & 0 & 0 & 0 &  &  & \\
    1 & 1 & 0 & 0 & 1 & 0 & 0 &  &  & \\
    1 & 1 & 0 & 0 & -1 & 0 & 0 &  &  & \\
    \hline
    & 1/2 & 1/2 & 0 & 1 & 0 & 0 & 0 & 1/2 & -1/2\\
  \end{array}}
\]
Метод \texttt{SRK1Wm} имеет сильный порядок
$(p_{d},p_{s}) = (1.0, 1.0)$, а метод \texttt{SRK2Wm} сильный порядок
$(p_{d},p_{s}) = (2.0, 1.0)$.


\subsection{Стохастические методы Рунге--Кутты слабой сходимости $p=2.0$}
Численные методы со слабой сходимостью хорошо аппроксимируют
характеристики распределения случайного процесса
$x^{\alpha}(t)$. Слабый численный метод не нуждается в информации о
траектории ведущего винеровского процесса $W^{\alpha}_{n}$ и случайные
величины для этих методов могут быть сгенерированы на другом
вероятностном пространстве. Поэтому можно воспользоваться легко
моделируемым на компьютере распределением.
\begin{align*}
& X^{0i\alpha} = x^{\alpha}_{n} + \sum\limits^{s}_{j=1}A_{0j}^{i}f^{\alpha}(t_{n} + c^{j}_{0}h_{n}, X^{0j\beta})h_{n} + \sum\limits^{s}_{j=1}\sum\limits^{m}_{l=1}B_{0j}^{i}G^{\alpha}_{l}(t_{n}+c^{j}_{1}h_{n}, X^{lj\beta})\hat{I}^{l},\\
& X^{ki\alpha} = x^{\alpha}_{n} + \sum\limits^{s}_{j=1}A^{i}_{1j}f^{\alpha}(t_{n}+c^{j}_{0}h_{n}, X^{0j\beta})h_{n} + \sum\limits^{s}_{j=1}B^{i}_{1j}G^{\alpha}_{k}(t_{n} + c_{1}^{j}h_{n}, X^{kj\beta})\sqrt{h_{n}},\\
& \widehat{X}^{ki\alpha} = x^{\alpha}_{n} + \sum\limits^{s}_{j=1}A^{i}_{2j}f^{\alpha}(t_{n}+c^{j}_{0}h_{n}, X^{0j\beta})h_{n} + \sum\limits^{s}_{j=1}\sum\limits^{m}_{l=1,l\neq k}B^{i}_{2j}G^{\alpha}_{l}(t_{n} + c_{1}^{j}h_{n}, X^{lj\beta})\frac{\hat{I}^{kl}}{\sqrt{h_{n}}},\\
& x^{\alpha}_{n+1} = x^{\alpha}_{n} + \sum\limits^{s}_{i=1}a_{i}f^{\alpha}(t_{n}+c^{i}_{0}, X^{0i\beta})h_{n} + \sum\limits^{s}_{i=1}\sum\limits^{m}_{k=1}\left(b^{1}_{i}\hat{I}^{k} + b^{2}_{i}\frac{\hat{I}^{kk}}{\sqrt{h_{n}}}\right) G^{\alpha}_{k}(t_{n}+c^{i}_{1}h_{n}, X^{ki\beta}) + \\
& +  \sum\limits^{s}_{i=1}\sum\limits^{m}_{k=1}\left(b^{3}_{i}\hat{I}^{k} + b^{4}_{i}\sqrt{h_{n}}\right)G^{\alpha}_{k}(t_{n}+c^{i}_{2}h_{n}, \widehat{X}^{ki\beta})
\end{align*}
Обобщенная таблица Бутчера имеет вид~\cite{L_Rossler_2010}:
\[
  {\renewcommand{\arraystretch}{1.5}%
  \begin{array}{c|c|c|c}
    c^{i}_{0} & A^{i}_{0j} & B^{i}_{0j} &\\ \hline
    c^{i}_{1} & A^{i}_{1j} & B^{i}_{1j} &\\ \hline
    c^{i}_{2} & A^{i}_{2j} & B^{i}_{2j} &\\ \hline
    & a_{i} & b^{1}_{i} & b^{2}_{i}\\ \hline
    &  & b^{3}_{i} & b^{4}_{i}\\
  \end{array}}
\]
{
\[
  {\renewcommand{\arraystretch}{1.5}%
  \begin{array}{c|ccc|ccc|ccc}
    0&0&0&0&0&0&0&&&\\
    1&1&0&0&\frac13&0&0&&&\\
    \frac{5}{12}&\frac{25}{144}&\frac{35}{144}&0&-\frac56&0&0&&&\\ \hline
    0&0&0&0&0&0&0&&&\\
    \frac14&\frac14&0&0&\frac12&0&0&&&\\
    \frac14&\frac14&0&0&-\frac12&0&0&&&\\ \hline
    0&0&0&0&0&0&0&&&\\
    0&0&0&0&1&0&0&&&\\
    0&0&0&0&-1&0&0&&&\\ \hline
    &\frac{1}{10}&\frac{3}{14}&\frac{24}{35}&1&-1&-1&0&1&-1\\ \hline
    &&&&\frac12&-\frac14&-\frac14&0&\frac12&-\frac12\\
  \end{array}}\quad
  {\renewcommand{\arraystretch}{1.5}%
  \begin{array}{c|ccc|ccc|ccc}
    0&0&0&0&0&0&0&&&\\
    1&1&0&0&1&0&0&&&\\
    0&0&0&0&0&0&0&&&\\ \hline
    0&0&0&0&0&0&0&&&\\
    1&1&0&0&1&0&0&&&\\
    1&1&0&0&-1&0&0&&&\\ \hline
    0&0&0&0&0&0&0&&&\\
    0&0&0&0&1&0&0&&&\\
    0&0&0&0&-1&0&0&&&\\ \hline
    &\frac{1}{2}&\frac{1}{2}&0&\frac12&\frac14&\frac14&0&\frac12&-\frac12\\ \hline
    &&&&-\frac12&\frac14&\frac14&0&\frac12&-\frac12\\
  \end{array}}
\]}

В численной схеме слабого стохастического метода Рунге-Кутты
используются следующие случайные величины:
\[
  \hat{I}^{kl} = 
  \left\{
    \begin{aligned}
      &\dfrac{1}{2}(\hat{I}^{k}\hat{I}^{l} - \sqrt{h}_{n}\tilde{I}^{k}),\; k<l,\\
      &\dfrac{1}{2}(\hat{I}^{k}\hat{I}^{l} + \sqrt{h}_{n}\tilde{I}^{l}),\; l<k,\\
      &\dfrac{1}{2}((\hat{I}^{k})^2 - h_{n}).\; k=l.
    \end{aligned}
  \right.
\]
Величина $\hat{I}^{k}$ имеет трехточечное распределение. Это означает,
что $\hat{I}^{k}$ может принимать три значения
$ \{-\sqrt{3h_{n}}, 0, \sqrt{3h_{n}}\}$ с вероятностями $1/6$, $2/3$ и
$1/6$ соответственно. В свою очередь величина $\tilde{I}^{k}$ имеет
двухточечное распределение $\{-\sqrt{h_{n}}, \sqrt{h_{n}}\}$ с
вероятностями $1/2$ и $1/2$.

\section{Анализ трудностей реализации стохастических численных методов Рунге--Кутты}

Как видно уже из формул, стохастические методы Рунге-Кутты значительно
сложнее своих классических аналогов. Кроме громоздкости формул, можно
выделить еще следующие факторы, усложняющие реализацию стохастических
методов в программном виде, а также их применение для численного
решения СДУ.
\begin{itemize}
\item При выборе конкретного метода надо учитывать какой тип
  сходимости необходимо обеспечить для данной конкретной задачи, а
  также какое из стохастических уравнений необходимо решать — в форме
  Ито или в форме Стратоновича. Это увеличивает количество алгоритмов,
  которые нужно запрограммировать.
\item Для методов с сильной сходимостью большей единицы на каждом шаге
  необходимо решать ресурсоемкую задачу по аппроксимации двукратных
  стохастических интегралов.
\item В численной схеме присутствуют не только матрицы и векторы, но и
  тензоры (четырехмерные массивы) с которыми необходимо совершать
  операцию свертки по нескольким индексам. Реализация свертки через
  суммирование с помощью обычных циклов приводит к существенному
  падению производительности.
\item Для использования слабых методов необходимо применять метод
  Монте Карло, проводя несколько серий по множеству испытаний в
  каждой. Так как метод Монте-Карло сходится приблизительно как
  $1/\sqrt{N}$, где $N$ --- число испытаний, то для достижения
  точности хотя бы $10^{-3}$, необходимо провести минимум $10^6$
  испытаний.
\end{itemize}

Наиболее существенное падение производительности происходит при
реализации универсального алгоритма, то есть такой программы, которая
может произвести расчет, используя произвольную таблицу
коэффициентов. В этом случае приходится использовать большое
количество вложенных циклов, для того, чтобы организовать
суммирование. Наличие в схемах для систем СДУ двойных сумм и сложной
комбинации индексов в множителях, находящихся под знаком этих сумм,
еще более усложняет задачу и число вложенных циклов вырастает до
шести.  Кроме этих специфических особенностей, стоит упомянуть также
несколько причин падения производительности, присущих также и
классическим схемам, которые в стохастическом случае также играют
роль. Очевидным способом хранения коэффициентов методов является
использование массивов. Однако у явных методов, которые мы
рассматриваем, матрица является нижне-диагональной и хранение ее в
виде двумерного массива приводит к тому, что более чем половина
выделенной для массива памяти тратится на хранение нулей.

Если изучить исходные коды популярных подпрограмм, реализующих
классические явные вложенные методы Рунге—Кутты, то можно обнаружит,
что в этих программах для хранения коэффициентов метода используется
набор именованных констант, а не массивов. Это вызвано также и тем,
что операции со скалярными величинами в большинстве языков
программирования проводятся быстрее, чем операции с массивами.

При сохранении требования универсальности создаваемого кода и вместе с
тем желание увеличить скорости вычислений и уменьшить расход памяти,
привели нас к решению использовать автоматическую генерацию кода по
одному шаблону для каждого отдельного метода.

Кроме выигрыша в производительности, автоматическая генерации кода,
позволяет добавлять или изменять все функции за раз путем
редактирования одного лишь шаблона, а не каждой отдельной функции. Это
позволяет как уменьшить количество ошибок, так и генерировать
различные варианты функций для разных целей.

\section{Краткое описание автоматической генерации кода}

В качестве языка, для написания генератора кода, используется язык
Python. Исходный код и документация созданных в нашем коллективе
программ, доступны по ссылке
\url{bitbucket.org/mngev/sde_num_generation}. Данный репозиторий
содержит модуль stochastic, в котором реализован винеровский случайный
процесс и численные методы типа Рунге-Кутты для скалярных и
многомерных стохастических дифференциальных уравнений в форме
Ито. Большая часть кода для этого модуля генерируется набором скриптов
из директории generator. В ручную написан класс, реализующий
многомерный винеровский процесс и методы Эйлера-Маруямы.

Для работы генератора кода мы использовали библиотеку для обработки
шаблонов Jinja2~\cite{L_Jinja2}. Данный шаблонизатор (template engine)
разрабатывался изначально для генерации HTML страниц, однако он
обладает очень гибким синтаксисом и может использоваться как
универсальное средство для генерации текстовых файлов любого вида, в
том числе и исходных кодов на любых языках программирования. Кроме
Jinja2 мы также использовали библиотеку NumPy для работы с массивами и
ускорения вычислений.

Кроме перечисленных выше двух внешних библиотек, был использован
стандартный модуль fraction, который позволяет задать коэффициенты
метода в виде рациональных дробей, а потом уже преобразовывать их в
вещественный вид с нужным порядком точности, а также модуль typing для
аннотации типов аргументов функций.

Шаблоны для генерации представляют собой файлы с исходным кодом на
языке Python со вставками специальных команд шаблонизатора
Jinja2. Информация о коэффициентах методов хранится отдельно, в
структурированном виде в формате JSON. Это позволяет легко добавлять
новые методы и изменять старые путем редактирования и JSON файлов. В
настоящее время генерируются методы с коэффициентами, представленными
в работах~\cite{L_Debrabant_2007,L_Debrabant_2013,L_Kloeden_Platen}.

В качестве языка для уже сгенерированных функций используется сам
python с активным применением библиотеки NumPy, которая позволяет
получить приемлемую производительность. Однако генерируемый код легко
можно переформатировать так, чтобы он соответствовал синтаксису любого
другого языка программирования. Мы планируем доработать программу, для
генерации кода на языке Julia (julialang.org). Данный язык был
представлен в 2012 году и изначально ориентирован на научные
вычисления. В настоящее время он интенсивно развивается и набирает
популярность. На сегодняшний день актуальной является версия
0.6.1. Julia позволяет получить производительность, сравнимую с C++ и
Fortran, но при этом является динамическим языком с возможностью
интерактивной командной строки (REPL) на подобии IPython и может
интегрироваться в интерактивную среду Jupyter. Текущая версия
библиотеки превосходит описанную авторами в
статье~\cite{kulyabov:2016:rk-stochastic}. Применение автогенерации позволило избавится от
большого количества вложенных циклов, выделения лишней памяти для
хранения нулевых коэффициентов, а также значительно упростило код.

\section{Заключение}
\label{sec:conclusion}

В данной работе рассмотрены стохастические численные схемы с порядком
сходимости выше 0.5. Показано, что такие методы существенно сложнее
эквивалентных численных методов для систем обыкновенных
дифференциальных уравнений. Выделены их особенности, делающие
эффективную программную реализацию таких методов не тривиальной
задачей. В заключительной части статьи обсуждается подход, основанный
на автоматической генерации программного кода, который позволяет
получить эффективную реализацию рассматриваемых методов, давая
возможность использовать любую таблицу коэффициентов. Дано краткое
описание программы, созданной авторами и ссылка на репозиторий с
исходными кодами.

\def\acknowledgmentsname{Благодарности}

\begin{acknowledgments}

Публикация подготовлена при поддержке Программы РУДН «5-100»
и при финансовой поддержке РФФИ в рамках научного проекта
\textnumero~16-07-00556.

\end{acknowledgments}

  \bibliography{rk-stoch,self}